\documentclass[11pt,a4paper]{article}
\PassOptionsToPackage{dvipsnames}{xcolor}
\usepackage[T1]{fontenc}
\usepackage{newtxtext}
\usepackage[margin=28mm,headheight=14pt]{geometry}
\usepackage{microtype}
\usepackage{amsmath,amsthm,mathtools}
\usepackage{newtxmath}
\usepackage{enumitem,booktabs,array}
\usepackage[dvipsnames]{xcolor}
\usepackage{fancyhdr}
\usepackage{aliascnt}
\usepackage{hyperref}
\usepackage[nameinlink,noabbrev]{cleveref}
\numberwithin{equation}{section}
\newtheorem{theorem}{Theorem}[section]
\newaliascnt{lemma}{theorem}
\newtheorem{lemma}[lemma]{Lemma}
\aliascntresetthe{lemma}
\newaliascnt{proposition}{theorem}

\aliascntresetthe{proposition}
\newaliascnt{corollary}{theorem}

\aliascntresetthe{corollary}
\theoremstyle{definition}
\newaliascnt{definition}{theorem}
\newtheorem{definition}[definition]{Definition}
\aliascntresetthe{definition}
\newaliascnt{remark}{theorem}
\newtheorem{remark}[remark]{Remark}
\aliascntresetthe{remark}
\crefname{lemma}{Lemma}{Lemmas}
\crefname{proposition}{Proposition}{Propositions}
\crefname{corollary}{Corollary}{Corollaries}
\crefname{definition}{Definition}{Definitions}
\crefname{remark}{Remark}{Remarks}
\crefname{theorem}{Theorem}{Theorems}
\newcommand{\R}{\mathbb R}

\newcommand{\E}{\mathcal E}
\newcommand{\F}{\mathcal F}
\newcommand{\D}{\mathcal D}
\newcommand{\Hrel}{\mathcal H}
\newcommand{\invN}{(-\Delta_N)^{-1}}

\newcommand{\dx}{\,\mathrm dx}
\newcommand{\dt}{\,\mathrm dt}

\newcommand{\norm}[1]{\left\|#1\right\|}

\newcommand{\mean}[1]{\overline{#1}}
\newcommand{\ceq}{c_*}

\allowdisplaybreaks[2]
\setlist{itemsep=2pt,topsep=5pt}

\title{\textbf {Global Dynamics of Finite-Energy Weak Solutions to a
Poisson--Nernst--Planck--Navier--Stokes System with
Interspecies Drag and Steric Effects}}
\author{Baoli Hao$^1$, Fanze Kong$^2$, and Chun Liu$^3$} 
\date{ }

\begin{document}

\maketitle
\begin{center}
$^1$ Department of Applied Mathematics, Illinois Institute of Technology, 10 West 35th Street,
Chicago, IL 60616 \\ {\tt bhao2@hawk.illinoistech.edu }

\noindent $^2$ Department of Applied Mathematics, University of Washington, 4182 West Stevens Way Northeast, Seattle, WA 98195\\
{\tt fzkong@uw.edu}

\noindent $^3$ Department of Mathematics,  City University of Hong Kong, Kowloon Tong, Hong Kong\\
{\tt Chun.Liu@cityu.edu.hk}

\end{center}

\renewcommand{\thefootnote}{\fnsymbol{footnote}}

\begin{abstract}
The Poisson--Nernst--Planck equations provide a continuum
framework for  ionic transport driven by diffusion
and electrostatic forces in electrolytes.
Using the energetic variational approach, we derive a
 Poisson--Nernst--Planck--Navier--Stokes
system that couples ionic transport with the flow of an
incompressible viscous solvent and incorporates steric
interactions and interspecies drag. The resulting ionic mobility matrix is non-diagonal, and the solvent dynamics are coupled to ionic transport through electrostatic forces. On a smooth bounded domain
$\Omega\subset\mathbb R^d$, $1\leq d\leq3$, we establish
the global existence of finite-energy weak solutions
for compatible initial data, including ionic vacuum.
We further establish exponential convergence of finite-energy weak solutions to the constant equilibrium using an entropy-dissipation inequality.
The existence proof uses a coupled implicit time
discretization, a Galerkin approximation of the solvent
velocity, and compactness arguments.
For  small perturbations of the
equilibrium, we also establish global existence and
uniqueness of strong solutions with exponential decay.
\end{abstract}
\noindent\textbf{Keywords.} Energetic Variational Approach; Interspecies Drag; Steric Effects; Navier-Stokes Flow Coupling; Global Weak Solutions; Entropy Dissipation.
 
\section{Introduction }
\label{sec:model}

Ionic transport is fundamental to many biological and
electrochemical systems, including cellular signaling \cite{Berridge2003},
energy storage \cite{TarasconArmand2001} and membrane separation 
\cite{Strathmann2010}.  The Poisson--Nernst--Planck (PNP) 
system provides a classical continuum description of ionic
 electrodiffusion by coupling the transport of charged species to
  the electrostatic potential \cite{Rubinstein1990,Eisenberg1998,Schuss2001}.
   However, classical PNP models typically rely on the dilute-solution 
   approximation, which neglects finite ion sizes and direct friction 
   between different ionic species. These effects can become important 
   in concentrated electrolytes and confined channels. 
   Steric interactions have been incorporated into PNP models through
    additional free-energy terms, which modify the chemical potentials 
    and introduce cross-diffusion
     effects \cite{Borukhov1997,Horng2012,LinEisenberg2014,LinLiuYiZhang2013}. 
     Interspecies friction provides a different coupling mechanism by
      penalizing the relative motion of 
      ionic species \cite{TaylorKrishna1993}. This effect can be
incorporated through an additional dissipation term,
leading to a non-diagonal, concentration-dependent
ionic mobility \cite{drag}.

{

The free-energy and dissipation structures also play
a central role in the mathematical analysis of these
PNP systems.  Entropy methods provide a general framework for 
establishing global weak solutions to cross-diffusion 
systems \cite{Jungel2015,Jungel2016}. Global weak existence 
results for PNP models with steric interactions have been established
 in \cite{Hsieh2019,GerstenmayerJungel2018,HirvonenJungel2024}.
  For models involving interspecies drag, Hsieh et al.~\cite{drag} 
  established local classical solvability. More recently, Hao, Kong, 
  Lam and Liu \cite{HKLL} proved the global existence of 
  finite-energy weak solutions in the presence of both steric
   interactions and interspecies drag, allowing ionic vacuum by 
   exploiting the coercivity of the steric free energy and weighted 
   gradient estimates arising from the dissipation.

Whereas, in many applications such as electro-osmosis and
electroconvection, ionic transport occurs within
a flowing solvent
\cite{SquiresBazant2004,Rubinstein2008}. The fluid advects the ions, and
the electric force on the charge density drives the fluid. This reciprocal  
interaction can be described by coupling the PNP equations with the 
incompressible Navier--Stokes equations \cite{Jerome2002,BotheFischerSaal2014}.  
The mathematical analysis of
classical PNP--NS systems is well developed.
In two dimensions, Ryham \cite{Ryham2009} established global weak
solvability, uniqueness and positive-time regularity for a binary system,
and Bothe, Fischer and Saal \cite{BotheFischerSaal2014} proved global strong
well-posedness and exponential convergence to equilibrium for multiple
species. Constantin and Ignatova \cite{ConstantinIgnatova} studied global
existence and long-time behavior under blocking and selective boundary
conditions, and Abdo and Ignatova \cite{AbdoIgnatova2024} obtained
exponential convergence in periodic domains. In three dimensions, Fischer
and Saal \cite{FischerSaal2017} obtained global weak solutions for $L^2$
data with species-dependent, nonconstant diffusion coefficients; Gong, Wang
and Zhang \cite{GongWangZhang2021} constructed suitable weak solutions and
proved partial regularity; and Constantin, Ignatova and Lee
\cite{CILNearEquilibrium} obtained global strong solutions near equilibrium.
Energy and entropy estimates also underlie the convergent fully discrete
schemes of Prohl and Schmuck \cite{ProhlSchmuck2010}.

}

{In this paper we couple the drag--steric PNP model of \cite{HKLL} to an
incompressible viscous solvent and derive the system and its
energy-dissipation law by the energetic variational approach
\cite{HyonKwakLiu2010}. Unlike \cite{HKLL}, the steric matrix is only
assumed positive semidefinite and coercive on the nonnegative cone. This
admits rank-one matrices, which arise when the steric energy depends only
on the total ionic volume fraction.
On smooth bounded domains $\Omega$ in dimensions
$1\le d\le3$, we construct global finite-energy weak solutions for
arbitrary compatible finite-energy data without any smallness assumption.
These solutions allow ionic vacuum. Both the concentrations and the fluid
velocity converge exponentially to the constant equilibrium. The
electrostatic force is square integrable in space and time for every
$d\le3$. In two dimensions this gives positive-time regularity of the
solvent velocity. Near equilibrium we also
obtain global existence, uniqueness and exponential decay of strong
solutions.}


\subsection{Model derivation via the energetic variational approach}
\label{sub:envara}

{We derive a coupled PNP--NS system incorporating steric interactions 
and interspecies drag via the energetic variational approach 
(EnVarA) \cite{HyonKwakLiu2010}. We first specify the total energy 
and dissipation. The total energy combines the solvent kinetic energy 
with the ionic free energy of \cite{HKLL}, comprising ideal mixing 
entropy, electrostatic energy, and local steric interactions. The 
dissipation accounts for ion--solvent friction, interspecies drag, 
and fluid viscosity. Then, the least action principle (LAP) yields 
the conservative forces, while the maximum dissipation principle (MDP) 
determines the dissipative forces from the dissipation. Finally, 
the force balances determine the ionic velocities and the fluid momentum 
equation. Substituting the ionic velocities into the mass conservation 
equations gives the Nernst--Planck equations. Together with the Poisson 
equation and the incompressibility condition, these equations form the 
coupled PNP--NS system. The detailed derivation is presented below.}

{
For each ionic species $i\in\{n,p\}$, let $c_i$ and $\mathbf v_i$
denote its concentration and transport velocity. The two species
carry charges $-z_0e$ and $+z_0e$, respectively, where $z_0>0$
is the common valence magnitude and $e>0$ the elementary charge. The solvent is an
incompressible Newtonian fluid with velocity $\mathbf u$. Thus,
\begin{equation}\label{eq:kinematics}
  \partial_t c_i+\nabla\cdot(c_i\mathbf v_i)=0,
  \qquad \nabla\cdot\mathbf u=0.
\end{equation}
The electrostatic potential satisfies
\begin{equation}\label{eq:dim-poisson}
  -\epsilon_e\Delta\phi=z_0e(c_p-c_n),
\end{equation}
where $\epsilon_e>0$ is the dielectric permittivity. For $f\in L^2(\Omega)$ with $\int_\Omega f\,\dx=0$, we write $\psi=\invN f$ for the unique solution of $-\Delta\psi=f$ in $\Omega$, $\partial_{\boldsymbol\nu}\psi=0$ on $\partial\Omega$, with $\int_\Omega\psi\,\dx=0$.

As shown in \cite{HKLL}, we model steric interactions by a quadratic
free-energy density $\frac12\mathbf c^\top G\mathbf c$, where
$\mathbf c=(c_n,c_p)^\top$ and
\(
G=\begin{pmatrix}g_{nn}&g_{np}\\g_{np}&g_{pp}\end{pmatrix}
\)
is symmetric positive semidefinite. Moreover, the Helmholtz free energy is
\begin{equation}\label{eq:free}
  \F(\mathbf c)=\int_\Omega\left[
    k_B\vartheta\sum_{i=n,p}c_i
      \left(\log\frac{c_i}{c_{\mathrm{ref}}}-1\right)
    +\frac12\mathbf c^\top G\mathbf c
    +\frac{\epsilon_e}{2}|\nabla\phi|^2
  \right]\dx,
\end{equation}
where $k_B$ is the Boltzmann constant, $\vartheta>0$ is the
absolute temperature, and $c_{\mathrm{ref}}>0$ is a reference concentration. The total energy
includes the solvent kinetic energy:
\begin{equation}\label{eq:total}
  \E(\mathbf c,\mathbf u)
  =\F(\mathbf c)+\frac{\rho_s}{2}\|\mathbf u\|_{L^2(\Omega)}^2,
\end{equation}
with solvent density $\rho_s>0$. Using the Poisson equation, the
electrochemical potentials are
\begin{equation}\label{eq:mu}
\begin{aligned}
  \mu_n&=k_B\vartheta\log\frac{c_n}{c_{\mathrm{ref}}}
            +g_{nn}c_n+g_{np}c_p-z_0e\phi,\\
  \mu_p&=k_B\vartheta\log\frac{c_p}{c_{\mathrm{ref}}}
            +g_{np}c_n+g_{pp}c_p+z_0e\phi.
\end{aligned}
\end{equation}
For a mass-preserving ionic displacement $\mathbf h_i$ with
$\mathbf h_i\cdot\boldsymbol\nu=0$ on $\partial\Omega$, 
\[
  \delta c_i=-\nabla\cdot(c_i\mathbf h_i),\qquad
  \delta_i\F=\int_\Omega c_i\nabla\mu_i\cdot\mathbf h_i\,\dx.
\]
LAP gives the conservative force density
\begin{equation}\label{eq:con-force}
  \mathbf F_i^{\mathrm{con}}=-c_i\nabla\mu_i.
\end{equation}
For the incompressible solvent, a divergence-free displacement
$\mathbf h_s$ vanishing on $\partial\Omega$ induces the velocity
variation $\delta\mathbf u=\partial_t\mathbf h_s
    +(\mathbf u\cdot\nabla)\mathbf h_s
    -(\mathbf h_s\cdot\nabla)\mathbf u.$
 Substituting this expression into the variation of the kinetic
action and integrating by parts, one formulates the inertial contribution as
\begin{align}\label{inertial-force-s}
  \mathbf F_s^{\mathrm{in}}
  =-\rho_s\bigl[\partial_t\mathbf u
  +(\mathbf u\cdot\nabla)\mathbf u\bigr].
\end{align}
The pressure enforces the incompressibility constraint.

For MDP, the dissipation rate as a function of the velocities is
\begin{equation}\label{eq:veldis}
\begin{aligned}
  \mathscr D(\mathbf c;\mathbf u,\mathbf v_n,\mathbf v_p)
  =\int_\Omega\Bigl[2\eta|\mathrm D\mathbf u|^2
    +\sum_{i=n,p}\frac{k_B\vartheta}{D_i}
      c_i|\mathbf v_i-\mathbf u|^2
    +\frac{k_B\vartheta}{D_{np}}c_nc_p
      |\mathbf v_n-\mathbf v_p|^2\Bigr]\dx,
\end{aligned}
\end{equation}
where $\mathrm D\mathbf u=(\nabla\mathbf u+\nabla\mathbf u^\top)/2$,
$\eta>0$ is the solvent viscosity, $D_n,D_p>0$ are ionic diffusion
coefficients, and $D_{np}>0$ characterizes interspecies drag.
Balancing the
resulting dissipative forces with the conservative and inertial
forces \eqref{eq:con-force}--\eqref{inertial-force-s} yield
\begin{align}
  \frac{k_B\vartheta}{D_n}c_n(\mathbf v_n-\mathbf u)
   +\frac{k_B\vartheta}{D_{np}}c_nc_p(\mathbf v_n-\mathbf v_p)
    &=-c_n\nabla\mu_n,\label{eq:force-n}\\
  \frac{k_B\vartheta}{D_p}c_p(\mathbf v_p-\mathbf u)
   +\frac{k_B\vartheta}{D_{np}}c_nc_p(\mathbf v_p-\mathbf v_n)
    &=-c_p\nabla\mu_p,\label{eq:force-p}\\
  \rho_s[\partial_t\mathbf u+(\mathbf u\cdot\nabla)\mathbf u]
       -\eta\Delta\mathbf u+\nabla p
    &=\sum_{i=n,p}\frac{k_B\vartheta}{D_i}
       c_i(\mathbf v_i-\mathbf u).\label{eq:solvent-raw}
\end{align}
 Defining the ionic fluxes relative to the solvent by
$\mathbf J_i=c_i(\mathbf v_i-\mathbf u)$, we solve \eqref{eq:force-n}--\eqref{eq:force-p}  to obtain
\begin{equation}\label{eq:flux}
  \begin{pmatrix}\mathbf J_n\\\mathbf J_p\end{pmatrix}
  =-M(\mathbf c)
    \begin{pmatrix}\nabla\mu_n\\\nabla\mu_p\end{pmatrix},
\end{equation}
where
\begin{align}\label{eq:M}
  &M(\mathbf c)=\frac{1}{k_B\vartheta\,\Xi(\mathbf c)}
  \begin{pmatrix}
    D_nc_n(D_{np}+D_pc_n)&D_nD_pc_nc_p\\
    D_nD_pc_nc_p&D_pc_p(D_{np}+D_nc_p)
  \end{pmatrix},\nonumber\\
  &\Xi(\mathbf c)=D_{np}+D_nc_p+D_pc_n.
\end{align}
Here the mobility matrix $M(\mathbf c)$ is symmetric
and positive semidefinite for $c_n,c_p\geq0$.
It is positive definite when $c_n,c_p>0$ and becomes
singular when either concentration vanishes.  Noting that 
\[
  \sum_{i=n,p}\frac{k_B\vartheta}{D_i}c_i(\mathbf v_i-\mathbf u)
  =-\nabla\left[k_B\vartheta(c_n+c_p)
      +\frac12\mathbf c^\top G\mathbf c\right]
   -z_0e(c_p-c_n)\nabla\phi,
\]
where we added the ionic force balances \eqref{eq:force-n}--\eqref{eq:force-p} and used \eqref{eq:mu}.
The gradient term can be absorbed into the pressure in
\eqref{eq:solvent-raw}.  Thus, combining the resulting momentum equation
with ionic conservation and \eqref{eq:dim-poisson}, we obtain
\begin{equation}\label{eq:system}
\left\{
\begin{aligned}
  \partial_t c_i+\mathbf u\cdot\nabla c_i
    &=\nabla\cdot\left(\sum_{j=n,p}M_{ij}(\mathbf c)\nabla\mu_j\right),
       &&i\in\{n,p\},\\
  -\epsilon_e\Delta\phi&=z_0e(c_p-c_n),\\
  \rho_s[\partial_t\mathbf u+(\mathbf u\cdot\nabla)\mathbf u]
       -\eta\Delta\mathbf u+\nabla p
    &=-z_0e(c_p-c_n)\nabla\phi,
       &&\nabla\cdot\mathbf u=0.
\end{aligned}
\right.
\end{equation}
Here $\mu_i$ and $M(\mathbf c)$ are given by \eqref{eq:mu}
and \eqref{eq:M}, respectively.

We next use the ionic force balances \eqref{eq:flux} to eliminate the
ionic velocities from the dissipation \eqref{eq:veldis} and express the
ionic contribution in terms of the concentrations
and chemical potential gradients. Noting that  the mobility matrix given by \eqref{eq:M} admits the decomposition
\begin{equation}\label{eq:decomposition}
  M(\mathbf c)
  =\frac{D_{np}}{k_B\vartheta\,\Xi(\mathbf c)}
    \operatorname{diag}(D_nc_n,D_pc_p)
   +\frac{D_nD_p}{k_B\vartheta\,\Xi(\mathbf c)}
    \mathbf c\otimes\mathbf c,
\end{equation}
 we have the ionic dissipation can be written as
\begin{equation}\label{eq:D}
\begin{aligned}
  \D(\mathbf c)
  &:=\int_\Omega \sum_{i,j=n,p}
    M_{ij}(\mathbf c)\nabla\mu_i\cdot\nabla\mu_j\,\dx\\
  &=\int_\Omega\frac{
    D_{np}D_nc_n|\nabla\mu_n|^2
    +D_{np}D_pc_p|\nabla\mu_p|^2
    +D_nD_p|c_n\nabla\mu_n+c_p\nabla\mu_p|^2
    }{k_B\vartheta\,\Xi(\mathbf c)}\,\dx.
\end{aligned}
\end{equation}
Substituting \eqref{eq:mu} into \eqref{eq:D}, we formally have from \eqref{eq:veldis} that  
\[
  \mathscr D(\mathbf c;\mathbf u,\mathbf v_n,\mathbf v_p)
  =2\eta\int_\Omega|\mathrm D\mathbf u|^2\,\dx+\D(\mathbf c).
\]
Thus, under the boundary conditions \eqref{eq:bc},
the system \eqref{eq:system} formally satisfies
\begin{equation}\label{eq:energyidentity}
  \frac{\mathrm d}{\mathrm dt}\E(\mathbf c,\mathbf u)
  +\mathscr D(\mathbf c;\mathbf u,\mathbf v_n,\mathbf v_p)=0.
\end{equation}
}
{
To nondimensionalize \eqref{eq:system}, we introduce
the Debye length and a reference diffusivity
\[
  \lambda_D
  =\left(
    \frac{\epsilon_e k_B\vartheta}
         {z_0^2e^2c_{\mathrm{ref}}}
  \right)^{1/2},
  \qquad
  D_*=\lambda_D
  \left(
    \frac{k_B\vartheta c_{\mathrm{ref}}}{\rho_s}
  \right)^{1/2}.
\]
We rescale the variables by
\[
\begin{gathered}
x=\lambda_D\widehat x,\quad
t=\frac{\lambda_D^2}{D_*}\widehat t,\quad
c_i=c_{\mathrm{ref}}\widehat c_i,\quad
\phi=\frac{k_B\vartheta}{z_0e}\widehat\phi, \quad
\mu_i=k_B\vartheta\widehat\mu_i,\quad
p=k_B\vartheta c_{\mathrm{ref}}\widehat p,
\end{gathered}
\]
and $(\mathbf u,\mathbf v_i,
\mathbf J_i/c_{\mathrm{ref}})
=\frac{D_*}{\lambda_D}
(\widehat{\mathbf u},\widehat{\mathbf v}_i,
\widehat{\mathbf J}_i)$.
With this choice of $D_*$, the coefficient of the fluid
inertial term in \eqref{eq:solvent-raw} is one, then we  define the dimensionless
parameters as follows:
\[
a=\frac{D_n}{D_*},\qquad
b=\frac{D_p}{D_*},\qquad
\delta=\frac{D_{np}}{D_*c_{\mathrm{ref}}},
\qquad
\nu=\frac{\eta}{\rho_sD_*}.
\]
Furthermore, we set the dimensionless steric matrix as
\begin{equation}\label{eq:Fmatrix}
F=\frac{c_{\mathrm{ref}}}{k_B\vartheta}G
=\begin{pmatrix}
f_{nn}&f_{np}\\
f_{np}&f_{pp}
\end{pmatrix}.
\end{equation}
Assume that $F$ is symmetric and positive semidefinite.  Moreover, there exists $\alpha>0$ such that
\begin{equation}\label{eq:F1}
  F=F^\top ,\qquad
  \mathbf r^\top F\mathbf r\geq\alpha|\mathbf r|^2
  \quad\text{for all }\mathbf r\in[0,\infty)^2.
\end{equation}
Hereafter, we omit the hats and denote the rescaled
domain by $\Omega$. In these variables, the ionic
mobility takes the form
\begin{equation}\label{eq:scaled-mobility}
M(\mathbf c)=
\frac{1}{\delta+ac_p+bc_n}
\begin{pmatrix}
ac_n(\delta+bc_n)&abc_nc_p\\
abc_nc_p&bc_p(\delta+ac_p)
\end{pmatrix}.
\end{equation}
Dividing the free and total energies by
$k_B\vartheta c_{\mathrm{ref}}\lambda_D^d$ gives
\begin{align}
\mathcal F(\mathbf c)
=\int_\Omega\left[
\sum_{i=n,p}c_i(\log c_i-1)
+\frac12\mathbf c^\top F\mathbf c
+\frac12|\nabla\phi|^2
\right]\,dx,
\quad
\mathcal E(\mathbf c,\mathbf u)
=\mathcal F(\mathbf c)
+\frac12\|\mathbf u\|_{L^2(\Omega)}^2.
\label{eq:scaled-total-energy}
\end{align}
The reduced ionic dissipation and the total dissipation are given by
\begin{equation}\label{nondimensional-dissipation-ionic}
\mathcal D(\mathbf c)
=\int_\Omega\sum_{i,j=n,p}
M_{ij}(\mathbf c)
\nabla\mu_i\cdot\nabla\mu_j\,dx, \qquad  \mathscr D(\mathbf c;\mathbf u,\mathbf v_n,\mathbf v_p)
  =2\nu\int_\Omega|\mathrm D\mathbf u|^2\,\dx+\D(\mathbf c).
\end{equation}
 Now, we have the following nondimensional system:
}
\begin{equation}\label{eq:nondimensional-system}
\left\{
\begin{aligned}
  &\partial_t c_i+\mathbf u\cdot\nabla c_i
  =\nabla\cdot\left(
    \sum_{j=n,p}M_{ij}(\mathbf c)\nabla\mu_j
  \right),
  \qquad i\in\{n,p\},\\[3pt]
  &-\Delta\phi=c_p-c_n,\\[3pt]
  &\partial_t\mathbf u
    +(\mathbf u\cdot\nabla)\mathbf u
    -\nu\Delta\mathbf u+\nabla p
    =-(c_p-c_n)\nabla\phi,\\[3pt]
  &\nabla\cdot\mathbf u=0,
\end{aligned}
\right.
\end{equation}
where
\begin{align}\label{nondimensional-pontential}
  \mu_n=\log c_n+f_{nn}c_n+f_{np}c_p-\phi,
  \qquad
  \mu_p=\log c_p+f_{np}c_n+f_{pp}c_p+\phi.
\end{align}
We impose the following boundary conditions
\begin{equation}\label{eq:bc}
 {\bf u}=0,\qquad \left(\sum_{j=n,p}M_{ij}({\bf c})\nabla\mu_j\right)
\cdot\boldsymbol{\nu}=0,\qquad \partial_{\boldsymbol{\nu}}\phi=0
 \quad\text{on }\partial\Omega,
 \qquad \int_\Omega\phi\,\dx=0.
\end{equation}
 The initial data satisfy
\begin{equation}\label{eq:data}
 c_{i,0}\geq0,\quad c_{i,0}\in L^2(\Omega),\quad
 \int_\Omega c_{n,0}\dx=\int_\Omega c_{p,0}\dx=m|\Omega|,
 \quad m>0,\qquad \mathbf u_0\in L^2_\sigma(\Omega),
\end{equation}
where
$L^2_\sigma(\Omega)=\{\mathbf v\in L^2(\Omega;\mathbb R^d):
\nabla\cdot\mathbf v=0 \text{ in }\Omega,\ \mathbf v\cdot\boldsymbol\nu=0
\text{ on }\partial\Omega\}$, with both conditions understood in the
weak sense; $L^r_\sigma(\Omega)$ is defined analogously for $1<r<\infty$.  
The homogeneous equilibrium is $(\mathbf c_*,0,0)$ with
\begin{align}\label{average-m}
  \mathbf c_*=(m,m)^\top,\qquad
  m=\frac{1}{|\Omega|}\int_\Omega c_{n,0}\,\dx
   =\frac{1}{|\Omega|}\int_\Omega c_{p,0}\,\dx .
\end{align}

The variational derivation preserves the non-diagonal ionic mobility and the
quadratic steric energy of \cite{HKLL}. Since the ideal and steric force
densities are gradients, they are absorbed into the pressure, so the
momentum equation in \eqref{eq:nondimensional-system} is driven only by the
electrostatic force.

The analysis of \eqref{eq:nondimensional-system} requires
compactness estimates that account for the non-diagonal
ionic mobility, steric cross-diffusion and fluid transport.
The first difficulty is that the energy bounds $\mathbf c$ only in
$L^\infty(0,T;L^2(\Omega)^2)$, while the dissipation controls the
concentration gradients only with weights. The degeneracy of $F$ makes
this harder than in \cite{HKLL}. When $F$ is positive definite, the matrix
$(\operatorname{diag}(c_n,c_p)^{-1}+F)^{-1}$ is bounded by $F^{-1}$ and the
individual dissipation suffices for the estimates of \cite{HKLL}. For
semidefinite $F$ this matrix can grow linearly in $c_n+c_p$.
Two complementary estimates resolve this. The individual fields $\mathbf A_n,\mathbf A_p$ control the
weighted concentration gradients
\begin{equation}\label{eq:intro-weighted}
 \int_0^T\!\!\int_\Omega
 \frac{|\nabla\mathbf c|^2}{(1+c_n+c_p)^2}\,
 \mathrm dx\,\mathrm dt
 \leq C_T.
\end{equation}
The collective field $\mathbf B$ controls the gradient of the ionic
pressure $c_n+c_p+\frac12\mathbf c^\top F\mathbf c$. This pressure is
comparable to $(1+c_n+c_p)^2$ on the nonnegative cone, so its bound gives
$\mathbf c\in L^4(Q_T)^2$. H\"older's inequality combines the two
estimates into $\mathbf c\in L^{4/3}(0,T;W^{1,4/3}(\Omega)^2)$, which gives
compactness of the concentrations. The collective dissipation is also
needed in the entropy-dissipation inequality that yields exponential decay. Poisson regularity also gives
$\nabla\phi\in L^\infty(0,T;L^4(\Omega;\mathbb R^d))$ for $d\le3$.
The electrostatic force $(c_p-c_n)\nabla\phi$ therefore lies in
$L^2(Q_T)$ in every dimension $d\le3$. In two dimensions this yields
positive-time regularity of the solvent velocity.

A second difficulty is to preserve the cancellation
between the contributions of ionic advection and
electrostatic forcing to the total energy balance.
We therefore employ a coupled implicit Euler--Galerkin
scheme in which both terms are evaluated at the same
time step. This cancellation yields uniform energy
and dissipation estimates for the discrete solutions.
Combined with bounds for the discrete time differences,
these estimates imply strong compactness of the
concentrations and velocities in
$L^2(0,T;L^2(\Omega)^2)$ \cite{Simon,DreherJungel}
and allow passage to the limit in $c_i\mathbf u$, $\mathbf u\otimes\mathbf u$,
and $(c_p-c_n)\nabla\phi$.

A further difficulty is the identification of the ionic
fluxes in the presence of vacuum. At each time step,
a mass-constrained entropy reconstruction yields strictly
positive concentrations with the prescribed ionic masses,
so that the logarithmic chemical potentials are well defined.
The positive lower bounds may depend on the approximation
parameters. We therefore identify the limiting fluxes
through weighted dissipation variables, obtaining flux
relations that remain valid where one or both concentrations
vanish.


\subsection{Main results}
\begin{definition}[]\label{def:weak}
Let $T>0$, write $Q_T=\Omega\times(0,T)$, and set
\begin{equation}\label{V-def-H-def}
 V=H^1_0(\Omega;\mathbb R^d)\cap L^2_\sigma(\Omega).
\end{equation}
A pair $(\mathbf c,\mathbf u)$, with
$\mathbf c=(c_n,c_p)^\top$, is called a finite-energy
weak solution to \eqref{eq:nondimensional-system} on $(0,T)$ if $c_n,c_p\ge0$ almost
everywhere in $Q_T$,
\begin{align*}
  c_n,c_p
  &\in L^\infty(0,T;L^2(\Omega))
    \cap L^{4/3}(0,T;W^{1,4/3}(\Omega)),\\
  \mathbf u
  &\in L^\infty(0,T;L^2_\sigma(\Omega))
    \cap L^2(0,T;V),
\end{align*}
and
\begin{equation}\label{eq:timeregs}
\begin{aligned}
  \partial_t c_n,\partial_t c_p
  \in L^{4/3}\bigl(0,T;(W^{1,4}(\Omega))'\bigr),\qquad \partial_t\mathbf u
  \in L^{4/3}(0,T;V').
\end{aligned}
\end{equation}
In addition, the species $c_n$ and $c_p$ satisfy
\[
  \int_\Omega c_i(t)\,\dx
  =\int_\Omega c_{i,0}\,\dx,
  \qquad i\in\{n,p\},
\]
for almost every $t\in(0,T)$, with the initial data
given in \eqref{eq:data}.
The electric potential is
$\phi=\invN(c_p-c_n)$, with homogeneous Neumann boundary
conditions and zero-mean normalization.

With all spatial gradients understood in the weak sense, let
\begin{equation}\label{eq:weights}
  \omega=\delta+a c_p+b c_n,\qquad
  \chi_n=\sqrt{\frac{\delta a c_n}{\omega}},\qquad
  \chi_p=\sqrt{\frac{\delta b c_p}{\omega}},\qquad
  \beta_i=c_i\sqrt{\frac{ab}{\omega}}
\end{equation}
and define
\begin{align}
\mathbf A_n
  &:=2\nabla\chi_n
    +\chi_n\nabla\bigl(
      \log\omega(\mathbf c)-\phi+f_{nn}c_n+f_{np}c_p
    \bigr),\notag\\
  \mathbf A_p
  &:=2\nabla\chi_p
    +\chi_p\nabla\bigl(
      \log\omega(\mathbf c)+\phi+f_{np}c_n+f_{pp}c_p
    \bigr),
    \label{eq:A}\\
  \mathbf B
  &:=\sqrt{\frac{b c_n}{\delta}}\,\mathbf A_n
    +\sqrt{\frac{a c_p}{\delta}}\,\mathbf A_p,
    \label{eq:B}\\
  \mathbf J_i
  &:=-\chi_i\mathbf A_i-\beta_i\mathbf B,
    \qquad i=n,p,
    \label{eq:Jweak}
\end{align}
 We require
 $ \chi_n,\chi_p
  \in L^{4/3}(0,T;W^{1,4/3}(\Omega))$,
$\mathbf A_n,\mathbf A_p\in L^2(Q_T;\mathbb R^d)$,
$\mathbf B\in L^2(Q_T;\mathbb R^d)$  
and 
$\mathbf J_i,c_i\mathbf u\in L^{4/3}(Q_T;\mathbb R^d)$.
Moreover, for every
$\psi_i\in L^4(0,T;W^{1,4}(\Omega))$, 
\begin{equation}\label{eq:weakions}
  \int_0^T
    \langle\partial_t c_i,\psi_i\rangle\,\dt
  -\int_{Q_T}
    (c_i\mathbf u+\mathbf J_i)\cdot\nabla\psi_i\,\dx\dt
  =0,
  \qquad i\in\{n,p\}.
\end{equation}
For every $\mathbf w\in L^4(0,T;V)$,
\begin{equation}\label{eq:weakns}
\begin{aligned}
  \int_0^T
    \langle\partial_t\mathbf u,\mathbf w\rangle\,\dt
   -\int_{Q_T}
    (\mathbf u\otimes\mathbf u):\nabla\mathbf w\,\dx\dt
   +\nu\int_{Q_T}
    \nabla\mathbf u:\nabla\mathbf w\,\dx\dt
   =-\int_{Q_T}
    (c_p-c_n)\nabla\phi\cdot\mathbf w\,\dx\dt.
\end{aligned}
\end{equation}
The initial data are attained in the sense that
\[
  c_i(t)\to c_{i,0}
  \quad\text{in }(W^{1,4}(\Omega))',
  \qquad i\in\{n,p\},
\]
and
\[
  \mathbf u(t)\to\mathbf u_0
  \quad\text{in }L^2_\sigma(\Omega),
  \qquad\text{as }t\downarrow0.
\]
Finally, for almost every $t\in(0,T)$, the total energy
satisfies
\begin{equation}\label{eq:weakenergy}
\begin{aligned}
  &\E(\mathbf c(t),\mathbf u(t))
   +\int_0^t\int_\Omega
    \left(
      |\mathbf A_n|^2+|\mathbf A_p|^2+|\mathbf B|^2
      +\nu|\nabla\mathbf u|^2
    \right)\,\dx\,ds
   \le \E(\mathbf c_0,\mathbf u_0).
\end{aligned}
\end{equation}
\end{definition}
We now state the existence of global finite-energy weak solutions.
\begin{theorem}
\label{thm:weak}
Let $\Omega\subset\mathbb R^d$ with $1\le d\le3$ be a smooth bounded domain,
and assume \eqref{eq:F1}, \eqref{eq:bc} and \eqref{eq:data}. Then there
exists a global solution $(\mathbf c,\mathbf u)$ whose restriction to every
$(0,T)$ is a finite-energy weak solution in the sense of
Definition~\ref{def:weak}. In addition, for every $T>0$,
\[
  \mathbf c\in L^4(Q_T)^2,\qquad
  \phi\in L^\infty(0,T;H^2(\Omega)),\qquad
  \mathbf J_i,\ c_i\mathbf u\in L^2(0,T;L^{4/3}(\Omega;\mathbb R^d)),
\]
and $\partial_tc_i\in L^2(0,T;(W^{1,4}(\Omega))')$. The electrostatic force
satisfies
\[
  (c_p-c_n)\nabla\phi\in L^\infty(0,T;L^{4/3}(\Omega;\mathbb R^d))\cap L^2(Q_T;\mathbb R^d).
\]
The initial data are attained strongly:
\[
  c_i(t)\to c_{i,0}
  \quad\text{in }L^2(\Omega),
  \qquad
  \mathbf u(t)\to\mathbf u_0
  \quad\text{in }L^2_\sigma(\Omega),
  \qquad t\downarrow0.
\]
Moreover, the energy inequality
\begin{equation}\label{eq:strongenergy}
  \E(\mathbf c(t),\mathbf u(t))
  +\int_s^t\!\!\int_\Omega
    \left(|\mathbf A_n|^2+|\mathbf A_p|^2+|\mathbf B|^2
    +\nu|\nabla\mathbf u|^2\right)\dx\,dr
  \le\E(\mathbf c(s),\mathbf u(s))
\end{equation}
holds for $s=0$ and almost every $s>0$ and for every $t\ge s$.
If $d=2$, then $\partial_t\mathbf u\in L^2(0,T;V')$ and for every $0<s<T$
\[
  \mathbf u\in C([s,T];V)\cap L^2(s,T;H^2(\Omega;\mathbb R^2)),
  \qquad
  \partial_t\mathbf u\in L^2(s,T;L^2_\sigma(\Omega)).
\]
If $\mathbf u_0\in V$, the same regularity holds with $s=0$.
\end{theorem}

\begin{remark}\label{rem:F}
Theorem~\ref{thm:weak} requires no smallness assumption on the initial data.
Condition~\eqref{eq:F1} includes all symmetric positive definite
matrices. It also admits rank-one matrices such as
$F=\left(\begin{smallmatrix}1&1\\1&1\end{smallmatrix}\right)$,
extending the positive-definiteness hypothesis in \cite{HKLL}.
\end{remark}

Our second main result shows that the constructed weak solutions converge
exponentially to the constant equilibrium.

\begin{theorem}
\label{thm:global-stability}
Let $(\mathbf c,\mathbf u)$ be a global weak solution constructed in
Theorem~\ref{thm:weak} and let $\phi=\invN(c_p-c_n)$. There exist $C,\kappa>0$
such that
\begin{equation}\label{eq:global-norm-decay}
\|{\bf c}(t)-{\bf c}_*\|_{L^2(\Omega)}
+\|{\bf u}(t)\|_{L^2(\Omega)}
+\|\phi(t)\|_{H^2(\Omega)}
\leq
C\bigl(\E(\mathbf c_0,\mathbf u_0)-\E(\mathbf c_*,0)\bigr)^{1/2}e^{-\kappa t/2},
\qquad t\geq0,
\end{equation}
where $\E(\mathbf c_0,\mathbf u_0)-\E(\mathbf c_*,0)\ge0$ is the initial
energy relative to the equilibrium. The constants $C$ and $\kappa$ depend
only on $\Omega,m,a,b,\delta,\nu,F,\alpha$ and an upper bound for
$\E(\mathbf c_0,\mathbf u_0)-\E(\mathbf c_*,0)$.
\end{theorem}


Our third main result addresses strong solutions near equilibrium.

\begin{theorem}\label{thm:small-strong}
Assume the hypotheses of Theorem~\ref{thm:weak} and fix $3<r<\infty$.
There exist $\varepsilon_*,C,\eta>0$, depending only on
$r,\Omega,m,a,b,\delta,\nu,F$, such that if
\[
  \mathbf c_0-\mathbf c_*\in W^{2-2/r,r}(\Omega)^2,
  \qquad
  \partial_{\boldsymbol\nu}\mathbf c_0=0,
  \qquad
  \mathbf u_0\in V\cap W^{2-2/r,r}(\Omega;\mathbb R^d),
\]
and
\begin{equation}\label{eq:small-strong-data}
  \|\mathbf c_0-\mathbf c_*\|_{W^{2-2/r,r}(\Omega)}
  +\|\mathbf u_0\|_{W^{2-2/r,r}(\Omega)}
  \leq\varepsilon_*,
\end{equation}
then \eqref{eq:nondimensional-system}--\eqref{eq:bc}
has a unique global strong solution satisfying, for every $T>0$,
\begin{equation}\label{eq:small-strong-regularity}
\begin{aligned}
  \mathbf c-\mathbf c_*
  &\in C([0,T];W^{2-2/r,r}(\Omega)^2)
     \cap W^{1,r}(0,T;L^r(\Omega)^2)
     \cap L^r(0,T;W^{2,r}(\Omega)^2),\\
  \mathbf u
  &\in C([0,T];W^{2-2/r,r}(\Omega;\mathbb R^d))
     \cap W^{1,r}(0,T;L^r_\sigma(\Omega))
     \cap L^r(0,T;W^{2,r}(\Omega;\mathbb R^d)).
\end{aligned}
\end{equation}
Moreover, $m/2\leq c_i(x,t)\leq3m/2$ for $i=n,p$, and
\begin{equation}\label{eq:small-strong-decay}
\begin{aligned}
  &\|\mathbf c(t)-\mathbf c_*\|_{W^{2-2/r,r}(\Omega)}
   +\|\mathbf u(t)\|_{W^{2-2/r,r}(\Omega)}
   +\|\phi(t)\|_{W^{4-2/r,r}(\Omega)}\\
  &\qquad\leq Ce^{-\eta t}
  \left(
    \|\mathbf c_0-\mathbf c_*\|_{W^{2-2/r,r}(\Omega)}
    +\|\mathbf u_0\|_{W^{2-2/r,r}(\Omega)}
  \right),
  \qquad t\geq0.
\end{aligned}
\end{equation}
\end{theorem}

The paper is organized as follows. Section~\ref{sec:prelim} collects
estimates for the Poisson equation and the Stokes--Galerkin projections and
proves weak lower semicontinuity of the energy. Section~\ref{sec:existence}
constructs global weak solutions through the coupled implicit
Euler--Galerkin scheme, identifies the limiting fluxes and proves the
additional two-dimensional fluid regularity. Section~\ref{sect4} proves an
entropy-dissipation inequality and the exponential convergence to
equilibrium. Section~\ref{sec:small-strong} treats strong solutions near
equilibrium, and Section~\ref{sect-final} discusses open problems.

\section{Preliminaries}
\label{sec:prelim}
This section presents the properties of the energy, the fluid estimates
and the a priori bounds used in Section~\ref{sec:existence}.
Set $\mathbf z=(-1,1)^\top$ and
\begin{align}\label{admissible-set-Km}
  \mathcal K_m:=\left\{
    \mathbf c\in L^2(\Omega)^2:
    c_i\geq0,\quad \int_\Omega c_i\,\dx=m|\Omega|,
    \quad i=n,p
  \right\}.
\end{align}
Thus, $c_p-c_n$ has zero mean for every $\mathbf c\in\mathcal K_m$.
\begin{lemma}
\label{lem:free-properties}
Let $\mathbf c\in\mathcal K_m$. Then $\F(\mathbf c)$ defined in \eqref{eq:scaled-total-energy} is finite and
\begin{equation}\label{eq:freecoercivityexpanded}
  \F(\mathbf c)\geq\frac\alpha2\|\mathbf c\|_2^2-2|\Omega|,
\end{equation}
where $\alpha$ is given in \eqref{eq:F1}.
Moreover, $\F$ is convex on $\mathcal K_m$. If
$\mathbf c_j\in\mathcal K_m$ and
$\mathbf c_j\rightharpoonup\mathbf c$ in $L^2(\Omega)^2$, then
$\mathbf c\in\mathcal K_m$ and
\begin{align}\label{lowersemicontinuity}
  \F(\mathbf c)\leq\liminf_{j\to\infty}\F(\mathbf c_j).
\end{align}
\end{lemma}
\begin{proof}
The bound follows from $r(\log r-1)\ge-1$ and \eqref{eq:F1}.
The mixing entropy is convex and grows at most quadratically. Both
quadratic terms are nonnegative and convex. Weak lower semicontinuity on the weakly closed set
$\mathcal K_m$ follows as in \cite[Lemmas~3.2--3.3]{HKLL}.
\end{proof}

We  next establish the  lower bound for the relative energy and show it is convex.
\begin{lemma}
\label{lem:relative-energy}
Let $\mathbf c\in\mathcal K_m$ and
$\mathbf u\in L^2(\Omega;\mathbb R^d)$. Set
\begin{align}\label{hm-def}
h_m(r):=r\log(r/m)-r+m,~~h_m(0)=m.
\end{align}
Then the relative energy defined by
\begin{equation}\label{H-energy}
\begin{aligned}
  \Hrel(\mathbf c,\mathbf u)
  &:=\E(\mathbf c,\mathbf u)-\E(\mathbf c_*,0)\\
  &=\frac12\|\mathbf u\|_{L^2(\Omega)}^2
   +\sum_{i=n,p}\int_\Omega h_m(c_i)\,\dx
   +\frac12\int_\Omega
     (\mathbf c-\mathbf c_*)^\top F(\mathbf c-\mathbf c_*)\,\dx
   +\frac12\|\nabla\phi\|_{L^2(\Omega)}^2.
\end{aligned}
\end{equation}
satisfies
\begin{equation}\label{eq:coerciveenergy}
  \Hrel(\mathbf c,\mathbf u)
  \geq\frac12\|\mathbf u\|_{L^2(\Omega)}^2
  +\gamma_*\|\mathbf c-\mathbf c_*\|_{L^2(\Omega)}^2
  +\frac12\|\nabla\phi\|_{L^2(\Omega)}^2,
\end{equation}
where $\gamma_*>0$ depends only on $m,F,\alpha$ given in \eqref{eq:F1}.  The functional $\Hrel$ is weakly lower semicontinuous on
$\mathcal K_m\times L^2(\Omega;\mathbb R^d)$.
Moreover, for $\mathbf c,\mathbf d\in\mathcal K_m$ and
$\mathbf u,\mathbf v\in L^2(\Omega;\mathbb R^d)$,
\begin{equation}\label{eq:midpointconvex}
\begin{aligned}
  \Hrel\left(\frac{\mathbf c+\mathbf d}{2},
             \frac{\mathbf u+\mathbf v}{2}\right)
  \leq{}&\frac{\Hrel(\mathbf c,\mathbf u)
                   +\Hrel(\mathbf d,\mathbf v)}2
  -\frac 18\|F^{1/2}(\mathbf c-\mathbf d)\|_{L^2(\Omega)}^2\\
  &-\frac18\|\mathbf u-\mathbf v\|_{L^2(\Omega)}^2-\frac14\sum_{i=n,p}\int_\Omega
\frac{(c_i-d_i)^2}{c_i+d_i}\,\mathrm dx.
\end{aligned}
\end{equation}
\end{lemma}

\begin{proof}
\eqref{H-energy} follows as in (4.1) of \cite{HKLL}
with the kinetic energy added.

The function $h_m$ is convex and attains its minimum
at $r=m$, where $h_m(m)=0$.
Together with $\log(r/m)\leq r/m-1$, this gives us
\begin{equation}\label{eq:hm-bounds}
0\leq h_m(r)\leq\frac{(r-m)^2}{m},
\qquad r\geq0,
\end{equation}
where the case $r=0$ follows by continuity.

To prove \eqref{eq:coerciveenergy}, let
$\mathbf r\in[0,\infty)^2$ with $|\mathbf r|$
sufficiently large. By \eqref{eq:F1}, we have
\[
\begin{aligned}
  \frac12(\mathbf r-\mathbf c_*)^\top
    F(\mathbf r-\mathbf c_*)
  &\geq\frac{\alpha}{2}|\mathbf r|^2
       -|F\mathbf c_*|\,|\mathbf r|\geq\frac{\alpha}{4}|\mathbf r|^2.
\end{aligned}
\]
On bounded sets of $\bf{r}$, we apply Taylor's formula to get from  $h_m''(r)=1/r$ that
  $\sum_i h_m(r_i)\geq C|\mathbf r-\mathbf c_*|^2$
for some $C>0$.
Combining these estimates with the positive semidefiniteness
of $F$ yields
\[
  \sum_{i=n,p}h_m(r_i)
  +\frac12(\mathbf r-\mathbf c_*)^\top
    F(\mathbf r-\mathbf c_*)
  \geq\gamma_*|\mathbf r-\mathbf c_*|^2
\]
for some $\gamma_*>0$.
Then, we integrate it to prove \eqref{eq:coerciveenergy}.

In light of Lemma \ref{lem:free-properties}, we have from the definition $\Hrel(\mathbf c,\mathbf u)
=\F(\mathbf c)-\F(\mathbf c_*)
+\frac12\|\mathbf u\|_{L^2(\Omega)}^2$ that $\Hrel$ is  weakly lower semicontinuous.   Finally, for fixed $t>0$, we perform the direct computation to obtain 
\[
  \frac{\partial^2}{\partial s^2}
  \left[
    \frac{h_m(s)+h_m(t)}2
    -h_m\left(\frac{s+t}{2}\right)
    -\frac{(s-t)^2}{4(s+t)}
  \right]
  =\frac{t(s-t)^2}{2s(s+t)^3}\geq0,
  \qquad s>0.
\]
Hence,
\begin{align}\label{elementary-identity-1}
  \frac{h_m(s)+h_m(t)}2
  -h_m\left(\frac{s+t}{2}\right)
  \geq\frac{(s-t)^2}{4(s+t)}
\end{align}
for all $s,t\geq0$.  Since $\mathbf c,\mathbf d\in\mathcal K_m$, we apply
\eqref{elementary-identity-1} with $s=c_i(x)$ and $t=d_i(x)$.
Then integrating over $\Omega$, summing over $i=n,p$, and using
\eqref{H-energy}, one obtains
\[
\begin{aligned}
  \Hrel\left(\frac{\mathbf c+\mathbf d}{2},
             \frac{\mathbf u+\mathbf v}{2}\right)
  \leq{}&
  \frac{\Hrel(\mathbf c,\mathbf u)
       +\Hrel(\mathbf d,\mathbf v)}2
   -\frac14\sum_{i=n,p}\int_\Omega
    \frac{(c_i-d_i)^2}{c_i+d_i}\,\dx\\
  &-\frac18\int_\Omega
    (\mathbf c-\mathbf d)^\top
    F(\mathbf c-\mathbf d)\,\dx
   -\frac18\|\mathbf u-\mathbf v\|_{L^2(\Omega)}^2.
\end{aligned}
\]
Since $F$ is symmetric and positive semidefinite, one finds
$(\mathbf c-\mathbf d)^\top F(\mathbf c-\mathbf d)
=|F^{1/2}(\mathbf c-\mathbf d)|^2$.
This proves \eqref{eq:midpointconvex}.


\end{proof}

The following inequality controls the change in free energy over one time step.
\begin{lemma}
\label{lem:bregman}
Let $\mathbf c,\mathbf d\in\mathcal K_m$ and assume that
$0<\underline c\leq c_i\leq\overline c<\infty$ almost everywhere
for $i=n,p$ with two constants $\underline c$ and $\overline c$. Define $\boldsymbol\mu(\mathbf c)
  :=\log\mathbf c+F\mathbf c+\mathbf z\,\invN(c_p-c_n).$
Then
\begin{equation}\label{eq:bregmanexact}
\begin{aligned}
  (\boldsymbol\mu(\mathbf c),\mathbf c-\mathbf d)
   -\bigl[\F(\mathbf c)-\F(\mathbf d)\bigr]&=\sum_{i=n,p}\int_\Omega
    \left[d_i\log\frac{d_i}{c_i}-d_i+c_i\right]\,\dx
   +\frac12\int_\Omega
    (\mathbf d-\mathbf c)^\top F(\mathbf d-\mathbf c)\,\dx\\
  &
   +\frac12\|\nabla \invN\bigl(\mathbf z\cdot(\mathbf d-\mathbf c)\bigr)\|_{L^2(\Omega)}^2
  \geq\frac 12\|F^{1/2}(\mathbf d-\mathbf c)\|_{L^2(\Omega)}^2.
\end{aligned}
\end{equation}

\end{lemma}

\begin{proof}
First of all, we apply the symmetry of $F$  to obtain 
\begin{align}\label{add1-lemma23}
\begin{aligned}
(F\mathbf c,\mathbf c-\mathbf d)
={}&\frac12(F\mathbf c,\mathbf c)
-\frac12(F\mathbf d,\mathbf d)+\frac12(F(\mathbf d-\mathbf c),\mathbf d-\mathbf c).
\end{aligned}
\end{align}
In addition, thanks to the self-adjointness of $\invN$ on mean-zero functions, one finds
\begin{align}\label{add2-lemma23}
\begin{aligned}
&\bigl(\mathbf z\cdot(\mathbf c-\mathbf d),
       \invN(\mathbf z\cdot\mathbf c)\bigr)\\
&\quad=
\frac12\|\nabla \invN(\mathbf z\cdot\mathbf c)\|_{L^2(\Omega)}^2
-\frac12\|\nabla \invN(\mathbf z\cdot\mathbf d)\|_{L^2(\Omega)}^2+
\frac12\|\nabla \invN(\mathbf z\cdot(\mathbf d-\mathbf c))\|_{L^2(\Omega)}^2.
\end{aligned}
\end{align}
Since
\[
(c_i-d_i)\log c_i
-c_i(\log c_i-1)+d_i(\log d_i-1)
=d_i\log\frac{d_i}{c_i}-d_i+c_i,
\qquad i=n,p,
\]
we integrate it  over $\Omega$ then sum it  over $i=n,p$,
and add \eqref{add1-lemma23} and \eqref{add2-lemma23}
to prove the equality in \eqref{eq:bregmanexact}.  Moreover, by using \eqref{eq:F1} and the nonnegativity of
the entropy and electrostatic terms, we obtain the
inequality in \eqref{eq:bregmanexact} holds.
 
\end{proof}

We are now ready to show the continuity of the energies.
\begin{lemma}\label{lem:energycontinuity}
Let $\mathbf c_j,\mathbf c\in\mathcal K_m$ and assume that
$\mathbf c_j\to\mathbf c$ strongly in $L^2(\Omega)^2$.
Then $\F(\mathbf c_j)\to\F(\mathbf c)$.
If also $\mathbf u_j\to\mathbf u$ strongly in
$L^2(\Omega;\mathbb R^d)$, then
$\Hrel(\mathbf c_j,\mathbf u_j)\to\Hrel(\mathbf c,\mathbf u)$.
\end{lemma}

\begin{proof}
    We follow the argument in
Proof of Theorem~4.4 in \cite{HKLL}, apply the 
strong $L^2$ convergence of ${\bf c}_j$ and \eqref{eq:hm-bounds} to obtain from 
 Vitali's theorem that
\[
  h_m(c_{i,j})\to h_m(c_i)
  \quad\text{in }L^1(\Omega),\qquad i=n,p,
\]
where $h_m$ is given in \eqref{hm-def}.  Since $\int_\Omega(c_{i,j}-c_i)\,\dx=0$, we further have
\[
\begin{aligned}
  \int_\Omega
  \bigl[c_{i,j}(\log c_{i,j}-1)
        -c_i(\log c_i-1)\bigr]\,\dx
  &=\int_\Omega
    \bigl[h_m(c_{i,j})-h_m(c_i)\bigr]\,\dx
    +\log m\int_\Omega(c_{i,j}-c_i)\,\dx\\
  &=\int_\Omega
    \bigl[h_m(c_{i,j})-h_m(c_i)\bigr]\,\dx
  \rightarrow0,
  \qquad i=n,p.
\end{aligned}
\]
Moreover, the strong $L^2$ convergence of ${\bf c}_j$ gives us
\[
\begin{aligned}
  \left|
    \frac12\int_\Omega\mathbf c_j^\top F\mathbf c_j\,\dx
    -\frac12\int_\Omega\mathbf c^\top F\mathbf c\,\dx
  \right|\leq
  C\|\mathbf c_j-\mathbf c\|_{L^2(\Omega)}
  \bigl(\|\mathbf c_j\|_{L^2(\Omega)}
        +\|\mathbf c\|_{L^2(\Omega)}\bigr)
  \rightarrow0,
\end{aligned}
\]
where constant $C>0$ depends only on $F$.
The standard elliptic regularity yields $\invN(\mathbf z\cdot\mathbf c_j)
  \to \invN(\mathbf z\cdot\mathbf c)$
\text{in }$H^2(\Omega)$,
which implies 
\[
  \frac12\|\nabla \invN(\mathbf z\cdot\mathbf c_j)\|_{L^2(\Omega)}^2
  \rightarrow
  \frac12\|\nabla \invN(\mathbf z\cdot\mathbf c)\|_{L^2(\Omega)}^2.
\]
Hence, $\F(\mathbf c_j)\to\F(\mathbf c)$ with $\F$ given by \eqref{eq:scaled-total-energy}.  Finally, strong $L^2$ convergence of $\mathbf u_j$ gives us
 $\frac12\|\mathbf u_j\|_{L^2(\Omega)}^2
  \rightarrow
  \frac12\|\mathbf u\|_{L^2(\Omega)}^2.$
Together with \eqref{H-energy}, this proves
$\Hrel(\mathbf c_j,\mathbf u_j)\to
\Hrel(\mathbf c,\mathbf u)$.
\end{proof}

Equip the space $V$ in \eqref{V-def-H-def} with $\|\mathbf v\|_V:=\|\nabla\mathbf v\|_{L^2(\Omega)}$.
For $d=2,3$, let $A$ be the Dirichlet Stokes operator associated with
$(\nabla\mathbf v,\nabla\mathbf w)$ on $V$. Its eigenfunctions
$(\mathbf e_r)_{r\ge1}$ form an orthonormal basis of $L^2_\sigma(\Omega)$,
with $A\mathbf e_r=\lambda_r^S\mathbf e_r$ and $\lambda_r^S>0$.
Let $P_N$ be the $L^2$-orthogonal projection onto
\begin{equation}\label{VN-approximation-space}
 V_N=\operatorname{span}\{\mathbf e_1,\ldots,\mathbf e_N\}.
\end{equation}
The eigenfunction expansion and Parseval's identity give
\begin{equation}\label{eq:projectionbounds}
  \|P_N\mathbf v\|_{L^2(\Omega)}\leq\|\mathbf v\|_{L^2(\Omega)}
  \quad(\mathbf v\in L^2_\sigma(\Omega)),\qquad
  \|P_N\mathbf v\|_V\leq\|\mathbf v\|_V
  \quad(\mathbf v\in V).
\end{equation}
In each space, $P_N\mathbf v\to\mathbf v$ in norm
\cite[Chapter~I, Section~2.6]{Temam}.
In dimension one, incompressibility and the no-slip condition give
$V=V_N=\{0\}$ and $P_N=0$.
For $\mathbf v_1,\mathbf v_2,\mathbf v_3\in V$, write
$b_{\rm NS}(\mathbf v_1,\mathbf v_2,\mathbf v_3)
=\int_\Omega(\mathbf v_1\cdot\nabla)\mathbf v_2\cdot\mathbf v_3\,\dx$.
Integration by parts gives antisymmetry in the last two arguments and
$b_{\rm NS}(\mathbf v_1,\mathbf v_2,\mathbf v_2)=0$.
For $\mathbf v,\mathbf w\in V$, interpolation \cite{Temam} and H\"older's inequality give
\begin{equation}\label{eq:ladyexpanded}
\|\mathbf v\|_{L^4(\Omega)}^2
\leq
\begin{cases}
C\|\mathbf v\|_{L^2(\Omega)}
 \|\nabla\mathbf v\|_{L^2(\Omega)},
& d=2,\\[2mm]
C\|\mathbf v\|_{L^2(\Omega)}^{1/2}
 \|\nabla\mathbf v\|_{L^2(\Omega)}^{3/2},
& d=3,
\end{cases}
  \qquad
  |b_{\rm NS}(\mathbf v,\mathbf v,\mathbf w)|
  \leq\|\mathbf v\|_{L^4(\Omega)}^2\|\nabla\mathbf w\|_{L^2(\Omega)}.
\end{equation}
Integrating in time and using Poincar\'e's inequality gives, for
$\mathbf u\in L^\infty(0,T;L^2_\sigma(\Omega))\cap L^2(0,T;V)$ and $1\le d\le3$,
\begin{align}\label{eq:fluid-L4-prelim}
\int_0^T \|\mathbf u(t)\|_{L^4(\Omega)}^{8/3}\,\mathrm dt
\leq
C\|\mathbf u\|_{L^\infty(0,T;L^2_\sigma(\Omega))}^{2/3}
 \|\mathbf u\|_{L^2(0,T;V)}^2.
\end{align}
If $d=2,$ we have
\begin{equation}\label{eq:fluid-L4-prelim-2d}
  \int_0^T\|\mathbf u(t)\|_{L^4(\Omega)}^4\,\dt
  \leq C\|\mathbf u\|_{L^\infty(0,T;L^2_\sigma(\Omega))}^2
         \|\mathbf u\|_{L^2(0,T;V)}^2.
\end{equation}

For positive concentrations the weighted fields have the following
expressions in terms of the chemical potentials.
\begin{lemma}
\label{lem:dissipation-representation}
Let $\mathbf c\in\mathcal K_m\cap W^{1,\infty}(\Omega)^2$
have a positive lower bound.   Then
\begin{align}\label{A-B-J-formula}
  \mathbf A_i=\chi_i\nabla\mu_i,\qquad
  \mathbf B=\sqrt{\frac{ab}{\omega}}
    \sum_{i=n,p}c_i\nabla\mu_i,\qquad
  \mathbf J_i=-\sum_{j=n,p}M_{ij}(\mathbf c)\nabla\mu_j,
\end{align}
where $\chi_i,\mathbf A_i,\mathbf B,\mathbf J_i$ are given by
\eqref{eq:weights}--\eqref{eq:Jweak}.
Moreover,
\begin{equation}\label{eq:Dfe}
  \D(\mathbf c)
  =\int_\Omega\left(
    |\mathbf A_n|^2+|\mathbf A_p|^2+|\mathbf B|^2
  \right)\,\dx,
\end{equation}
where $\D$ is defined in \eqref{nondimensional-dissipation-ionic}. 
\end{lemma}

\begin{proof}
By using \eqref{eq:weights}, we differentiate to get 
$2\nabla\chi_i=\chi_i(\nabla\log c_i-\nabla\log\omega)$.
Substituting it into \eqref{eq:A}, we have
$\mathbf A_i=\chi_i\nabla\mu_i$.  It then follows from \eqref{eq:B}  that 
the formula for $\mathbf B$ in \eqref{A-B-J-formula} holds.  In light of \eqref{eq:scaled-mobility}, we decompose
$M(\mathbf c)
  =\operatorname{diag}(\chi_n^2,\chi_p^2)
    +\frac{ab}{\omega}\mathbf c\otimes\mathbf c$, then obtain \eqref{A-B-J-formula} 
 and \eqref{eq:Dfe} hold.
\end{proof}

\begin{lemma}
\label{lem:coercivity}
For $\mathbf c\in(0,\infty)^2$, let
$H(\mathbf c):=\operatorname{diag}(c_n^{-1},c_p^{-1})+F$.
Then there exists $\kappa>0$, depending only on $a,b,\delta,F$,
such that
\begin{equation}\label{eq:HMH}
  H(\mathbf c)M(\mathbf c)H(\mathbf c)
  \geq\frac{\kappa}{(1+c_n+c_p)^2}I.
\end{equation}
Let $\mathbf c\in\mathcal K_m\cap W^{1,\infty}(\Omega)^2$.  Assume $\operatorname*{ess\,inf}_\Omega c_i>0$. Set
$\phi=\invN(c_p-c_n)$ and
$\boldsymbol\mu=\log\mathbf c+F\mathbf c+\mathbf z\phi$.
If $\Hrel(\mathbf c,0)\leq \Hrel(\mathbf c_0,\mathbf u_0)$ with $\Hrel$ given in \eqref{H-energy}, then
\begin{equation}\label{eq:weighted}
\begin{aligned}
  I_w(\mathbf c)
  :=\int_\Omega\frac{|\nabla\mathbf c|^2}{(1+c_n+c_p)^2}\,\dx
  \leq C\D(\mathbf c)
       +C\int_\Omega |\nabla\phi|^2\,\dx
   \leq C\D(\mathbf c)+C_{\Hrel(\mathbf c_0,\mathbf u_0)},
\end{aligned}
\end{equation}
where
$\D(\mathbf c)$ is given in \eqref{eq:Dfe} and $C,C_{\Hrel(\mathbf c_0,\mathbf u_0)}>0$ are constants independent of the positive lower bounds of $c_i$.
\end{lemma}
\begin{proof}

Since $F$ is positive semidefinite, we expand to obtain 
\[
\begin{aligned}
  H(\mathbf c)\operatorname{diag}(c_n,c_p)H(\mathbf c)
  &=
  \operatorname{diag}(c_n^{-1},c_p^{-1})
  +2F+F\operatorname{diag}(c_n,c_p)F\geq\frac1{c_n+c_p}I.
\end{aligned}
\]
By using \eqref{eq:weights}, one further has
\[
  H(\mathbf c)\operatorname{diag}(\chi_n^2,\chi_p^2)
  H(\mathbf c)
  \geq\frac{\delta\min\{a,b\}}{\omega(c_n+c_p)}I
  \geq\frac{\kappa}{(1+c_n+c_p)^2}I.
\]
It then follows from  $M(\mathbf c)\geq\operatorname{diag}(\chi_n^2,\chi_p^2)$ that
 \eqref{eq:HMH} holds.
Moreover,
$H(\mathbf c)\geq\operatorname{diag}(c_n^{-1},c_p^{-1})$
and \eqref{eq:HMH} imply
\begin{equation}\label{eq:cone-inverse}
\begin{aligned}
  \|H(\mathbf c)^{-1}\|_{\mathrm{op}}
  &\leq c_n+c_p,\\
  \left\|
    H(\mathbf c)^{-1}
    \operatorname{diag}(\chi_n^{-1},\chi_p^{-1})
  \right\|_{\mathrm{op}}
  &\leq C(1+c_n+c_p).
\end{aligned}
\end{equation}
Invoking \eqref{A-B-J-formula}, we get 
$\nabla\boldsymbol\mu
=H(\mathbf c)\nabla\mathbf c+\mathbf z(\nabla\phi)^\top$.  Then, we use \eqref{eq:cone-inverse} to get
\begin{equation}\label{eq:cone-gradient-pointwise}
  \frac{|\nabla\mathbf c|}{1+c_n+c_p}
  \leq C\bigl(
    |\mathbf A_n|+|\mathbf A_p|+|\nabla\phi|
  \bigr),
\end{equation}
which indicates $I_w(\mathbf c)
  \leq C\D(\mathbf c)+C\|\nabla\phi\|_{L^2(\Omega)}^2,$
where we used \eqref{eq:Dfe}.
Finally, combining \eqref{H-energy} with
$\|\nabla\phi\|_{L^2(\Omega)}^2\leq2\Hrel(\mathbf c_0,\mathbf u_0)$, we 
completes the proof of \eqref{eq:weighted}.


\end{proof}

The flux bound below remains valid when one species vanishes.
\begin{lemma}
\label{lem:vacuum-flux}
Let $T>0$ and $\mathbf c$ be nonnegative 
with $\mathbf c(t)\in\mathcal K_m$ for almost every $t$.
Assume
\[
  \mathbf c\in L^\infty(0,T;L^2(\Omega)^2)
    \cap L^{4/3}(0,T;W^{1,4/3}(\Omega)^2),
  \qquad
  \chi_i\in L^{4/3}(0,T;W^{1,4/3}(\Omega)),
\]
and that the fields in \eqref{eq:A}--\eqref{eq:B} satisfy
$\mathbf A_n,\mathbf A_p,\mathbf B\in L^2(Q_T;\mathbb R^d)$.
Then, $\mathbf J_i$  given in \eqref{eq:Jweak} satisfies
$\mathbf J_i\in L^2(0,T;L^{4/3}(\Omega;\mathbb R^d))$ and
\begin{equation}\label{eq:weak-flux-bound}
\begin{aligned}
  \|\mathbf J_i\|_{L^2(0,T;L^{4/3}(\Omega))}
  \leq C\left(
    \|\mathbf A_i\|_{L^2(Q_T)}
    +\|c_i\|_{L^\infty(0,T;L^2(\Omega))}^{1/2}
      \|\mathbf B\|_{L^2(Q_T)}
  \right).
\end{aligned}
\end{equation}
Furthermore, $\mathbf A_i=\mathbf J_i=0$ almost everywhere
on $\{c_i=0\}$.
Thus, $\D(\mathbf c)$ in \eqref{eq:Dfe} is well-defined.
\end{lemma}

\begin{proof}
Thanks to \eqref{eq:weights}, we have
\begin{equation}\label{eq:weightbounds}
\begin{gathered}
  0\leq\chi_n\leq\sqrt{\frac{\delta a}{b}},
\quad
  0\leq\chi_p\leq\sqrt{\frac{\delta b}{a}},\quad
  0\leq\beta_n\leq\sqrt{ac_n},
  \quad
  0\leq\beta_p\leq\sqrt{bc_p}.
\end{gathered}
\end{equation}
Similarly as the argument shown in the proof of  Lemma~3.7 in \cite{HKLL}, we obtain
\[
  \|\mathbf J_i(t)\|_{L^{4/3}(\Omega)}
  \leq C\left(
    \|\mathbf A_i(t)\|_{L^2(\Omega)}
    +\|c_i(t)\|_{L^2(\Omega)}^{1/2}
     \|\mathbf B(t)\|_{L^2(\Omega)}
  \right).
\]
Taking the $L^2$ norm in time then proves
\eqref{eq:weak-flux-bound}.  Moreover, as shown in Remark~3.1 of \cite{HKLL}, one finds
$\chi_i=\nabla\chi_i=0$ almost everywhere on
$\{c_i=0\}$, so \eqref{eq:A} gives us $\mathbf A_i=0$
there. Since $\beta_i=0$ on the same set, we use 
\eqref{eq:Jweak} to get $\mathbf J_i=0$.
Finally, the assumed $L^2(Q_T)$ bounds on
$\mathbf A_n,\mathbf A_p,\mathbf B$ ensure that
$\D(\mathbf c)$ defined by \eqref{eq:Dfe}
belongs to $L^1(0,T)$.
\end{proof}

To control the electric force we need more than the energy bound on
$\mathbf c$ in $L^\infty(0,T;L^2)$. The collective field $\mathbf B$
controls the gradient of the sum of the ideal and steric pressures and yields the additional
integrability of $\mathbf c$ used below.
\begin{lemma}
\label{lem:mixed}
Let $T>0$ and $(\mathbf c,\mathbf u)$ be measurable functions
such that, for almost every $t\in(0,T)$,
\[
  \mathbf c(t)\in\mathcal K_m\cap W^{1,\infty}(\Omega)^2,
  \qquad \operatorname*{ess\,inf}_\Omega c_i(t)>0,
  \qquad \mathbf u(t)\in V.
\]
Let $\phi=\invN(c_p-c_n)$ and define $\chi_i,\mathbf A_i,\mathbf B,
\mathbf J_i$ by \eqref{eq:weights}--\eqref{eq:Jweak}.
Assume
\begin{equation}\label{eq:aprioribase}
  \operatorname*{ess\,sup}_{0<t<T}\Hrel(\mathbf c(t),\mathbf u(t))
  +\int_0^T\left(
    \D(\mathbf c(t))+\nu\|\nabla\mathbf u(t)\|_2^2
  \right)\,\dt\leq C_0
\end{equation}
for some constant $C_0>0,$ where $\D$ is given in \eqref{eq:Dfe}.  Then
\begin{align}
  &\mathbf c\in L^\infty(0,T;L^2(\Omega))
      \cap L^{4/3}(0,T;W^{1,4/3}(\Omega)),\qquad
    \phi\in L^\infty(0,T;H^2(\Omega)),
    \label{eq:c-reg}\\
  &\chi_i\in L^\infty(Q_T)\cap L^{4/3}(0,T;W^{1,4/3}(\Omega)),\qquad
    \mathbf A_n,\mathbf A_p,\mathbf B\in L^2(Q_T),
    \label{eq:chi-reg}\\
  &\mathbf J_i,\ c_i\mathbf u\in L^2(0,T;L^{4/3}(\Omega)),\qquad
    \mathbf u\in L^\infty(0,T;L^2_\sigma(\Omega))\cap L^2(0,T;V),
    \label{eq:flux-reg}\\
  &\mathbf c\in L^4(Q_T)^2,\qquad
    (c_p-c_n)\nabla\phi\in L^\infty(0,T;L^{4/3}(\Omega;\R^d))\cap L^2(Q_T;\R^d).
    \label{eq:force-reg}
\end{align}
The norms in \eqref{eq:c-reg}--\eqref{eq:force-reg}
are bounded by a constant $C=C(T,C_0,\Omega,a,b,\delta,\nu,F,m)$,
independent of the positive lower bounds and the
$W^{1,\infty}$ norms of ${\bf c}$.
\end{lemma}
\begin{proof}
By using \eqref{eq:aprioribase} and \eqref{eq:coerciveenergy} in Lemma \ref{lem:relative-energy}, one finds constant $C>0$ such that
\begin{align}\label{c-L2x-Linf-t-lem210}
  \|\mathbf c(t)\|_{L^2(\Omega)}
  \leq \left(\frac{2C_0}{\gamma_*}\right)^{1/2}
       +m\sqrt{2|\Omega|}
  \leq C
\end{align}
for almost every $t\in(0,T)$.  Standard elliptic regularity then gives us
\begin{align}\label{H2-standard-elliptic-regularity-C}
  \|\phi(t)\|_{H^2(\Omega)}
  \leq C\|c_p(t)-c_n(t)\|_{L^2(\Omega)}
  \leq C.
\end{align}
for some constant $C>0.$  Moreover,  
\eqref{eq:weighted} in Lemma \ref{lem:coercivity} and \eqref{eq:aprioribase} yield for a constant $C>0,$
\begin{align}\label{weighted-estimate-Iw-C}
  \int_0^T I_w(\mathbf c(t))\,\dt
  \leq C\int_0^T\D(\mathbf c(t))\,\dt+TC_{C_0}
  \leq C,
\end{align}
where $C_{C_0}>0$ is a constant depending on $C_0.$

By \eqref{eq:F1}, the total ionic pressure satisfies
\begin{align}\label{upper-lower-bound-cone}
1+c_n+c_p+\frac12\mathbf c^\top F\mathbf c
  \asymp(1+c_n+c_p)^2.
  \end{align}
Its gradient is determined by the collective field in \eqref{A-B-J-formula}:
$\nabla\left(c_n+c_p+\frac12\mathbf c^\top F\mathbf c\right)
  =\sqrt{\frac{\omega}{ab}}\,\mathbf B
   -(c_p-c_n)\nabla\phi.$
Consequently, \eqref{eq:weights} and \eqref{upper-lower-bound-cone} give
\begin{align}\label{grad-energy-cn-cp-F}
  \left|
    \nabla\left(
      1+c_n+c_p+\frac12\mathbf c^\top F\mathbf c
    \right)^{3/4}
  \right|^2
  \leq C|\mathbf B|^2
       +C(c_n+c_p)|\nabla\phi|^2,
\end{align}
where $C>0$ is  a constant.  Then, \eqref{c-L2x-Linf-t-lem210} and \eqref{H2-standard-elliptic-regularity-C} imply for a constant $C>0,$
\begin{align}\label{Holder-weighted-gradphi}
\int_\Omega(c_n+c_p)|\nabla\phi|^2\,\dx
  \leq\|c_n+c_p\|_{L^2(\Omega)}
       \|\nabla\phi\|_{L^4(\Omega)}^2
  \leq C.
\end{align}
By collecting   \eqref{grad-energy-cn-cp-F}, \eqref{eq:aprioribase} and \eqref{Holder-weighted-gradphi}, one has
\[
\begin{aligned}
  \int_0^T\!\!\int_\Omega
  \left|
    \nabla\left(
      1+c_n+c_p+\frac12\mathbf c^\top F\mathbf c
    \right)^{3/4}
  \right|^2\,\dx\,\dt
  &\leq
  C\int_{Q_T}|\mathbf B|^2\,\dx\,\dt
  +C\int_{Q_T}
    (c_n+c_p)|\nabla\phi|^2\,\dx\,\dt\\
  &\leq CC_0+CT.
\end{aligned}
\]
Moreover, in light of \eqref{c-L2x-Linf-t-lem210} and \eqref{eq:data}, one finds
\[
  \operatorname*{ess\,sup}_{0<t<T}
  \int_\Omega
  \left(1+c_n+c_p+\frac12\mathbf c^\top F\mathbf c\right)
  \,\dx\leq C.
\]
Combining these estimates, we invoke Poincar\'e's inequality to obtain
\[
  \int_0^T
  \left\|
    \left(1+c_n+c_p+\frac12\mathbf c^\top F\mathbf c\right)^{3/4}
  \right\|_{H^1(\Omega)}^2\,\dt
  \leq C,
\]
where $C>0$ is  a constant.  Since $H^1(\Omega)\hookrightarrow L^4(\Omega)$,
we apply the interpolation inequality to further get  for constant $ C>0,$
\[
\begin{aligned}
  &\int_{Q_T}
  \left(1+c_n+c_p+\frac12\mathbf c^\top F\mathbf c\right)^2
  \,\dx\,\dt\\
  &\leq
  \left[
    \operatorname*{ess\,sup}_{0<t<T}
    \int_\Omega
    \left(1+c_n+c_p+\frac12\mathbf c^\top F\mathbf c\right)
    \,\dx
  \right]^{1/2}
  \int_0^T
  \left\|
    \left(1+c_n+c_p+\frac12\mathbf c^\top F\mathbf c\right)^{3/4}
  \right\|_{L^4(\Omega)}^2\,\dt\\
  &\leq C,
\end{aligned}
\]
Consequently, \eqref{eq:F1} implies that for constant $C>0,$
\begin{equation}\label{eq:cone-higher-integrability}
  \|c_n\|_{L^4(Q_T)}+\|c_p\|_{L^4(Q_T)}\leq C.
\end{equation}
By using H\"older's inequality, we obtain from
\eqref{weighted-estimate-Iw-C}
and \eqref{eq:cone-higher-integrability} that
\begin{equation}\label{eq:mixedholder}
\begin{aligned}
  \|\nabla c_i\|_{L^{4/3}(Q_T)}
  &\leq
  \left\|\frac{\nabla c_i}{1+c_n+c_p}\right\|_{L^2(Q_T)}
  \|1+c_n+c_p\|_{L^4(Q_T)}
  \leq C,
  \qquad i=n,p.
\end{aligned}
\end{equation}
Together with \eqref{c-L2x-Linf-t-lem210}, this proves \eqref{eq:c-reg}. 

In light of \eqref{eq:Dfe} and \eqref{eq:aprioribase}, we obtain
$\mathbf A_n,\mathbf A_p,\mathbf B$
are bounded in $L^2(Q_T)$.
Since $\omega=\delta+bc_n+ac_p\geq\delta$,
\eqref{eq:A} and \eqref{eq:weightbounds} imply
\[
  |\chi_i|\leq C,\qquad
  |\nabla\chi_i|
  \leq C\bigl(
    |\mathbf A_i|+|\nabla\mathbf c|+|\nabla\phi|
  \bigr),
  \qquad i=n,p,
\]
for constant $C>0$.  Thus, we have  \eqref{eq:chi-reg} holds.
Furthermore, we use  \eqref{eq:weak-flux-bound} in
Lemma~\ref{lem:vacuum-flux} to get
\[
  \|\mathbf J_i\|_{L^2(0,T;L^{4/3}(\Omega))}
  \leq C,
  \qquad i=n,p,
\]
where $C>0$ is a constant.  We next estimate the velocity and the coupling terms.   Thanks to  \eqref{eq:coerciveenergy} in Lemma \ref{lem:relative-energy} and \eqref{eq:aprioribase}, one finds
\[
  \frac12\|\mathbf u\|_{L^\infty(0,T;L^2_{\sigma}(\Omega))}^2
  +\nu\int_0^T
    \|\nabla\mathbf u(t)\|_{L^2(\Omega)}^2\,\dt
  \leq C_0.
\]
The embedding $V\hookrightarrow L^4(\Omega;\mathbb R^d)$ gives
\[
  \|\mathbf u\|_{L^\infty(0,T;L^2_{\sigma}(\Omega))}
  +\|\mathbf u\|_{L^2(0,T;V)}
  +\|\mathbf u\|_{L^2(0,T;L^4(\Omega))}
  \leq C
\]
for constant $C>0.$
Combining this with \eqref{eq:c-reg}, we obtain
\[
  \|c_i\mathbf u\|_{L^2(0,T;L^{4/3}(\Omega))}
  \leq
  \|c_i\|_{L^\infty(0,T;L^2(\Omega))}
  \|\mathbf u\|_{L^2(0,T;L^4(\Omega))}
  \leq C,
\]
where $C>0$ is a constant.  This proves \eqref{eq:flux-reg}.

Estimate~\eqref{eq:cone-higher-integrability} is the concentration bound in
\eqref{eq:force-reg}. It also gives $\mathbf c\in L^2(0,T;L^4(\Omega)^2)$
since $T<\infty$. For $d\le3$, the embedding
$H^2(\Omega)\hookrightarrow W^{1,4}(\Omega)$ and \eqref{eq:c-reg}
give $\nabla\phi\in L^\infty(0,T;L^4(\Omega;\mathbb R^d))$.
Hence, by H\"older's inequality,
\[
\begin{aligned}
  &\|(c_p-c_n)\nabla\phi\|
    _{L^\infty(0,T;L^{4/3}(\Omega;\mathbb R^d))}
\leq
  \|c_p-c_n\|_{L^\infty(0,T;L^2(\Omega))}
  \|\nabla\phi\|_{L^\infty(0,T;L^4(\Omega;\mathbb R^d))}
  \leq C,
\end{aligned}
\]
where $C>0$ is a constant. The same bounds yield
\begin{equation}\label{eq:forceL2}
  \|(c_p-c_n)\nabla\phi\|_{L^2(Q_T)}
  \leq
  \|c_p-c_n\|_{L^2(0,T;L^4(\Omega))}
  \|\nabla\phi\|_{L^\infty(0,T;L^4(\Omega))}
  \leq C,
\end{equation}
which completes the proof of \eqref{eq:force-reg}.
\end{proof}

\section{Global finite-energy weak solutions}
\label{sec:existence}
In this section we prove Theorem~\ref{thm:weak}. The construction combines
an implicit Euler discretization in time, a Galerkin approximation of the
velocity and a regularization in entropy variables. At each time step a
constrained minimization recovers positive concentrations with the
prescribed masses, and the discrete problem is solved by a fixed-point
argument. Uniform energy and dissipation bounds give compactness, and the
limiting fluxes are identified through the weighted dissipation variables.
We extend the approximation in Section~3.3 of \cite{HKLL} by solving
the ionic equations and the Galerkin momentum equation at the same
time step. The ionic unknowns are entropy variables. To recover positive
concentrations with the prescribed masses from a bounded, componentwise
mean-zero field $\mathbf y$, set
\begin{equation}\label{eq:reconstructmin}
 { c}[{\bf y}]=\underset{{\bf c}\in\mathcal K_m}{\operatorname{argmin}}
        \left\{\F({\bf c})-\int_{\Omega} {\bf y}\cdot {\bf c}\, \dx\right\},\qquad
 \phi[{\bf y}]=\invN({\bf z}\cdot { c}[{\bf y}]),
\end{equation}
Here $\mathcal K_m$ enforces the two mass constraints. The reconstruction
inverts the local entropy gradient.
\begin{lemma}\label{lem:localinverse}
Define $h({\bf r})=\sum_i r_i(\log r_i-1)+{\bf r}^\top F{\bf r}/2$ with $F$ given in \eqref{eq:Fmatrix}, we have the map $\nabla h:(0,\infty)^2\to\R^2$ is a bijection with a smooth inverse ${\bf g}$. Moreover,
\begin{equation}\label{eq:gderivative}
 D{\bf g}({\boldsymbol\xi})=\bigl[D^2h({\bf g}(\boldsymbol\xi))\bigr]^{-1},\qquad
 |{\bf g}(\boldsymbol\xi)-{\bf g}(\boldsymbol\eta)|\leq L_B|\boldsymbol{\xi}-\boldsymbol{\eta}|,~|\xi|,|\eta|\leq B,
\end{equation}
where $L_B:=
\sup_{|\boldsymbol\xi|\leq B}
\left\|
\bigl[D^2h(g(\boldsymbol\xi))\bigr]^{-1}
\right\|
<\infty.$  Moreover, the convex conjugate $h^*$ is smooth, with
$\nabla h^*={\bf g}$ and a positive definite Hessian.
\end{lemma}

\begin{proof}
 For any $\boldsymbol\xi\in\mathbb R^2$, we consider the function ${\bf r}\mapsto h({\bf r})-\boldsymbol\xi\cdot{\bf r}$.  Thanks to \eqref{eq:F1}, we have the function 
is coercive and strictly convex on $[0,\infty)^2$. It therefore has a unique minimizer.  We next show that  the minimizer lies in $(0,\infty)^2$.
Indeed, if its $i$th component vanishes, then
\[
h({\bf r}+\varepsilon e_i)
-\boldsymbol\xi\cdot({\bf r}+\varepsilon e_i)
-\bigl(h({\bf r})-\boldsymbol\xi\cdot{\bf r}\bigr)
=\varepsilon\log\varepsilon+O(\varepsilon)<0
\]
for sufficiently small $\varepsilon>0$, which is  contradicted to
minimality.  Thus,  the minimizer satisfies
$\nabla h({\bf r})=\boldsymbol\xi$.  In addition, the strict convexity of ${\bf r}\mapsto h({\bf r})-\boldsymbol\xi\cdot{\bf r}$ ensures that $\nabla h({\bf r})=\boldsymbol\xi$ has at most
one solution, so $\nabla h$ is a bijection.

We next prove \eqref{eq:gderivative}.  Thanks to \eqref{eq:F1}, we obtain that $D^2h({\bf r})
=\operatorname{diag}(r_n^{-1},r_p^{-1})+F
$ is positive definite.  Then, we apply 
the inverse function theorem to get that the inverse of $\nabla h$ denoted by ${\bf g}$
is smooth.  Moreover, differentiating
$\nabla h({\bf g}(\boldsymbol\xi))=\boldsymbol\xi$ gives us $D{\bf g}({\boldsymbol\xi})=\bigl[D^2h({\bf g}(\boldsymbol{\xi}))\bigr]^{-1}$.
Then we integrate it along the line segment joining
$\boldsymbol\xi$ and $\boldsymbol\eta$ to conclude \eqref{eq:gderivative}.

Finally, by using the definition of the convex conjugate, one has  $h^*(\boldsymbol\xi)
=\boldsymbol\xi\cdot {\bf g}(\boldsymbol\xi)
-h({\bf g}(\boldsymbol\xi)).$
We then differentiate to obtain $\nabla h^*(\boldsymbol\xi)={\bf g}(\boldsymbol\xi)$ and  the Hessian  is positive definite.
\end{proof}

\begin{lemma}\label{lem:reconstruct}
The minimizer in problem \eqref{eq:reconstructmin} exists and is unique. Moreover, there is a unique constant vector $\boldsymbol{\zeta}$ such that
\begin{equation}\label{eq:reconstruction}
 \log c[{\bf y}]+Fc[{\bf y}]+{\bf z}\phi[{\bf y}]={\bf y}+{\boldsymbol \zeta}.
\end{equation}
Also, for $\norm {\bf y}_{L^\infty(\Omega)}\leq R$ with constant $R>0$, there are positive constants $\underline c_R$, $\overline c_R$ and  $C_R$ independent of ${\bf y}$ such that 
\begin{equation}\label{eq:reconstructionbounds}
 0<\underline c_R\leq c_i[{\bf y}]\leq\overline c_R<\infty,
 \quad |\boldsymbol{\zeta}|+\norm{\phi[{\bf y}]}_{H^2(\Omega)}\leq C_R.
\end{equation} 
Furthermore, the reconstruction  shown in \eqref{eq:reconstructmin} is continuous from $L^\infty_0(\Omega)^2$ into $L^\infty(\Omega)^2\times H^2(\Omega)\times\R^2$, where $L^\infty_0(\Omega)$ denotes the space of 
bounded functions with zero  mean. If $y\in H^\ell(\Omega)$ with $\ell>1+\frac{d}{2}$, then $c[{\bf y}]\in W^{1,\infty}(\Omega)$.
\end{lemma}
\begin{proof}
\emph{Step 1. Existence and uniqueness.}
The admissible set $\mathcal K_m$ is nonempty, convex and weakly closed
in $L^2(\Omega)^2$. The objective in \eqref{eq:reconstructmin} is coercive
by \eqref{eq:freecoercivityexpanded} and weakly lower semicontinuous by
Lemma~\ref{lem:free-properties}. It therefore attains its minimum. The
strict convexity of the mixing entropy gives uniqueness. Since $\mathbf y$
has zero mean, comparison with $\mathbf c_*$ gives
\[
 \F(\mathbf c[\mathbf y])-\int_\Omega \mathbf y\cdot\mathbf c[\mathbf y]\,\dx
 \le \F(\mathbf c_*).
\]
For $\|\mathbf y\|_\infty\le R$, coercivity and Cauchy--Schwarz yield
$\|\mathbf c[\mathbf y]\|_2\le C_R$. Neumann regularity then bounds
$\phi[\mathbf y]$ in $H^2(\Omega)\hookrightarrow L^\infty(\Omega)$.

\emph{Step 2. Reduction to a local problem.}
Fix $\mathbf c=\mathbf c[\mathbf y]$ and $\phi=\phi[\mathbf y]$. Define
$L_{\mathbf y}(\mathbf d)=\int_\Omega[h(\mathbf d)
+(\mathbf z\phi-\mathbf y)\cdot\mathbf d]\,\dx$.
Expanding the electrostatic energy gives, for $\mathbf d\in\mathcal K_m$,
\begin{equation}\label{eq:localreduction}
 \F({\bf d})-\int_{\Omega} {\bf {y}}\cdot {\bf d}\,\dx-\left(\F({\bf c})-\int_{\Omega} {\bf y}\cdot {\bf c}\, \dx\right)
 =L_{{\bf y}}({\bf d})-L_{{\bf y}}({\bf c})+\tfrac12\norm{\nabla \invN({\bf z}\cdot({\bf d}-{\bf c}))}_{L^2(\Omega)}^2,
\end{equation}
Suppose $L_{\mathbf y}(\mathbf d)<L_{\mathbf y}(\mathbf c)$ and take
$\mathbf d_\theta=(1-\theta)\mathbf c+\theta\mathbf d$. By convexity,
$L_{\mathbf y}(\mathbf d_\theta)-L_{\mathbf y}(\mathbf c)
\le\theta[L_{\mathbf y}(\mathbf d)-L_{\mathbf y}(\mathbf c)]$.
The last term in \eqref{eq:localreduction} is $O(\theta^2)$ at
$\mathbf d_\theta$. For small $\theta>0$ the objective would therefore
decrease, contradicting minimality. Thus $\mathbf c$ minimizes $L_{\mathbf y}$.

\emph{Step 3. Mass multipliers.}
The two multipliers must enforce the prescribed masses.
Let ${\bf q}(x)={\bf y}(x)-{\bf z}\phi(x)$, a bounded function, and define
\begin{align}\label{Def-G-zeta}
G(\boldsymbol{\zeta})=\int_{\Omega} h^*({\bf q}+\boldsymbol{\zeta})\dx-m|\Omega|(\boldsymbol{\zeta}_n+\boldsymbol{\zeta}_p),
\end{align}
where $h^*$ is the convex conjugate of $h$ given in
Lemma~\ref{lem:localinverse}. By the same lemma, we have the inverse ${\bf g}$ is smooth
and $\nabla h^*={\bf g}$.  Moreover,
\begin{align}\label{G-derivative}
 \partial_{\zeta_i}G=\int_{\Omega} g_i({\bf q}+\boldsymbol{\zeta})\dx-m|\Omega|,
 \qquad D^2G=\int_{\Omega} D^2h^*({\bf q}+\boldsymbol{\zeta})\dx>0.
\end{align}
We prove that a minimizer of problem $\min_{\boldsymbol\zeta\in\mathbb R^2}G(\boldsymbol\zeta)
$ exists.  To this end, we fix $0<r_*<m$ and use the trial vector
${\bf d}_\zeta={\bf \ceq}+r_*\boldsymbol{\zeta}/|\boldsymbol{\zeta}|$ for $\boldsymbol{\zeta}\ne0$. By the definition of the convex conjugate $h^*({\bf q}+\boldsymbol\zeta)
=\sup_{{\bf r}\in[0,\infty)^2}
\bigl\{({\bf q}+\boldsymbol\zeta)\cdot{\bf r}-h({\bf r})\bigr\}$, we take ${\bf r}={\bf d}_{\boldsymbol\zeta}$ to get
\[
h^*({\bf q}+\boldsymbol\zeta)
\geq({\bf q}+\boldsymbol\zeta)\cdot{\bf d}_{\boldsymbol\zeta}
-h({\bf d}_{\boldsymbol\zeta}).
\]
Since $m-r_*\leq({\bf d}_{\boldsymbol\zeta})_i\leq m+r_*$ for 
 $i=n,p$, we further have from \eqref{Def-G-zeta} that
\[
 G(\boldsymbol{\zeta})\geq\int_{\Omega}[({\bf q}+\boldsymbol{\zeta})\cdot \boldsymbol{d}_\zeta-h(\boldsymbol{d}_\zeta)]\,\dx
                -|\Omega|\boldsymbol{\zeta}\cdot{\bf \ceq}
 \geq r_*|\Omega||\boldsymbol{\zeta}|-C_R,
\]
where $C_R>0$ is a constant.  Since $G$ given in \eqref{Def-G-zeta} is continuous,
 we have the minimizer of $\min_{\boldsymbol\zeta\in\mathbb R^2}G(\boldsymbol\zeta)
$ exists.
Moreover, the positive definiteness of $D^2G$ shown in \eqref{G-derivative}
implies that $G$ is strictly convex, so this
minimizer is unique.

Let $\boldsymbol\zeta$ be the minimizer of $\min_{\boldsymbol\zeta\in\mathbb R^2}G(\boldsymbol\zeta)
$.  By using \eqref{G-derivative},  the first-order optimality conditions yield
\[
0=\partial_{\zeta_i}G(\boldsymbol\zeta)
=\int_\Omega g_i({\bf q}+\boldsymbol\zeta)\,\dx-m|\Omega|,
\qquad i=n,p.
\]
Therefore, $\widehat{\bf c}:=g({\bf q}+\boldsymbol\zeta)$
satisfies the prescribed mass constraints.
Moreover, thanks to Lemma \ref{lem:localinverse}, we have $\widehat{\bf c}$ is bounded and strictly
positive.  
Thus $\widehat{\bf c}\in\mathcal K_m$.

It remains to identify $\widehat{\bf c}$ with ${\bf c}$.
For every ${\bf d}\in\mathcal K_m$, we apply Fenchel's inequality  to get
\[
h({\bf d})-{\bf q}\cdot{\bf d}
\geq
\boldsymbol\zeta\cdot{\bf d}
-h^*({\bf q}+\boldsymbol\zeta),
\]
where the
equality holds when ${\bf d}=\widehat{\bf c}$ due to
$\nabla h(\widehat{\bf c})={\bf q}+\boldsymbol\zeta$.
By integrating and using $\int_\Omega d_i\,\dx
=\int_\Omega\widehat c_i\,\dx
=m|\Omega|$,
one finds from the definition of $L_{\bf y}$ that
\[
\begin{aligned}
L_{\bf y}({\bf d})
\geq
\boldsymbol\zeta\cdot\int_\Omega{\bf d}\,\dx
-\int_\Omega h^*({\bf q}+\boldsymbol\zeta)\,\dx
=
\boldsymbol\zeta\cdot\int_\Omega\widehat{\bf c}\,\dx
-\int_\Omega h^*({\bf q}+\boldsymbol\zeta)\,\dx=L_{\bf y}(\widehat{\bf c}).
\end{aligned}
\]
By Step 2 and strict convexity, $\widehat{\mathbf c}$ and $\mathbf c$
are the same minimizer of $L_{\mathbf y}$. This proves \eqref{eq:reconstruction}.
For fixed $\mathbf y$, the argument $\mathbf q+\boldsymbol\zeta$ of
$\mathbf g$ is bounded. Continuity of $\mathbf g$ bounds $\mathbf c$,
and \eqref{eq:reconstruction} then bounds $\log c_i$ for $i=n,p$.
 
\emph{Step 4. Uniform bounds.}
Assume $\|\mathbf y\|_\infty\le R$. Averaging \eqref{eq:reconstruction}
and applying Jensen's inequality gives
\[
\zeta_i=\mean{\log c_i}+(F{\bf\ceq})_i\leq\log m+(F{\bf \ceq})_i,
\]
where the overbar denotes spatial averaging over $\Omega$, componentwise for vector fields. It gives an upper bound independent of ${\bf y}$ for $\boldsymbol{\zeta}$. Taking the scalar product of \eqref{eq:reconstruction} with ${\bf c}(x)$,  since ${\bf y}-{\bf z}\phi$ is bounded by $C_R$, we have
\begin{align}\label{C-upper-lower-bound}
 \alpha|{\bf c}|^2+\sum_i c_i\log c_i\leq C_R(c_n+c_p)\leq \sqrt{2}C_R|{\bf c}|
 \end{align}
for constant $C_R>0.$  The lower bound $r\log r\geq-e^{-1}$ then yields a uniform upper bound for $|{\bf c}|$.  For $i\in\{n,p\}$, since $|q_i|$ is bounded, we have from \eqref{C-upper-lower-bound} that $|q_i-(F{\bf c})_i|\leq C_R$ for some constant $C_R>0.$   Moreover, thanks to \eqref{eq:reconstruction}, we obtain 
\begin{align}\label{ci-average-bound}
  e^{\zeta_i-C_R}\leq c_i(x)\leq e^{\zeta_i+C_R}
  \qquad\text{for almost every }x\in\Omega.
\end{align}
Taking spatial averages in \eqref{ci-average-bound} and using
$|\Omega|^{-1}\int_\Omega c_i\,dx=m$, we  get $\log m-C_R\leq\zeta_i\leq\log m+C_R.$
It then follows that
\[
  m e^{-2C_R}\leq c_i(x)\leq m e^{2C_R}
  \qquad\text{for almost every }x\in\Omega.
\]
Together with the $H^2$ bound for $\phi$ shown in Step 1, this proves
\eqref{eq:reconstructionbounds}.  

\emph{Step 5. Continuity.}
For inputs $\mathbf y_1,\mathbf y_2$, the two reconstructed concentrations
have the same masses. Hence $\int_\Omega
  \bigl(\boldsymbol{\zeta}[{\bf y}_1]-\boldsymbol{\zeta}[{\bf y}_2]\bigr)
  \cdot\bigl(c[{\bf y}_1]-c[{\bf y}_2]\bigr)\,\dx=0.$
Moreover, since $\phi$ satisfies the Poisson equation, we have
\[
  \int_\Omega
  \bigl(\phi[{\bf y}_1]-\phi[{\bf y}_2]\bigr)
  {\bf z}\cdot\bigl(c[{\bf y}_1]-c[{\bf y}_2]\bigr)\,\dx
  =
  \norm{\nabla(\phi[{\bf y}_1]-\phi[{\bf y}_2])}_{L^2(\Omega)}^2.
\]
By using \eqref{eq:reconstruction} and these identities, we obtain
\begin{align*}
  &\sum_{i\in\{n,p\}}\int_\Omega
  \bigl(\log c_i[{\bf y}_1]-\log c_i[{\bf y}_2]\bigr)
  \bigl(c_i[{\bf y}_1]-c_i[{\bf y}_2]\bigr)\,\dx\\
  &\quad+\int_\Omega
  \bigl(c[{\bf y}_1]-c[{\bf y}_2]\bigr)^\top
  F\bigl(c[{\bf y}_1]-c[{\bf y}_{2}]\bigr)\,\dx\\
  &\quad+\norm{\nabla(\phi[{\bf y}_1]-\phi[{\bf y}_2])}_{L^2(\Omega)}^2\\
  &=\int_\Omega
  ({\bf y}_1-{\bf y}_2)\cdot\bigl(c[{\bf y}_1]-c[{\bf y}_2]\bigr)\,\dx.
\end{align*}
Since $\log$ is increasing, one further obtains from \eqref{eq:reconstructionbounds} and \eqref{eq:F1} that
\begin{equation}\label{eq:reconstructionlip}
 {\overline c}^{-1}_R\norm{c[{\bf y}_1]-c[{\bf y}_2]}_{L^2(\Omega)}^2
  +\norm{\nabla(\phi[{\bf y}_1]-\phi[{\bf y}_2])}_{L^2(\Omega)}^2
  \leq
  \int_\Omega
  ({\bf y}_1-{\bf y}_2)\cdot\bigl(c[{\bf y}_1]-c[{\bf y}_2]\bigr)\,\dx.
\end{equation}

Now suppose ${\bf y}_j\to {\bf y}$ in $L^\infty(\Omega)^2$.
By using \eqref{eq:reconstructionlip} and the
Cauchy--Schwarz inequality, we have
\[
  \norm{c[{\bf y}_j]-c[{\bf y}]}_{L^2(\Omega)}
  \leq {\overline c}_R\norm{{\bf y}_j-{\bf y}}_{L^2(\Omega)}
  \rightarrow0
\]
as $j\rightarrow+\infty.$
Invoking the elliptic regularity and the embedding
$H^2(\Omega)\hookrightarrow L^\infty(\Omega)$, we further obtain 
\[
  \norm{\phi[{\bf y}_j]-\phi[{\bf y}]}_{H^2(\Omega)}
  +\norm{\phi[{\bf y}_j]-\phi[{\bf y}]}_{L^\infty(\Omega)}
  \leq C\norm{c[{\bf y}_j]-c[{\bf y}]}_{L^2(\Omega)}
  \rightarrow0,
\]
where $C>0$ is a constant.
Since $(y_j)$ is bounded in $L^\infty(\Omega)^2$, we use 
\eqref{eq:reconstructionbounds} to show that there is a constant
$\underline c_R>0$ such that
\[
  c_i[{\bf y}_j]\geq\underline c_R,
  \qquad c_i[{\bf y}]\geq\underline c_R
  \qquad(i\in\{n,p\}).
\]
Thus, integrating \eqref{eq:reconstruction} over $\Omega$
and dividing by $|\Omega|$, we obtain
\[
  \zeta_i[{\bf y}_j]-\zeta_i[{\bf y}]
  =\mean{\log c_i[{\bf y}_j]-\log c_i[{\bf y}]}.
\]
 Since $|\log a-\log b|
  \leq\frac{|a-b|}{\underline c_R}$ with $a,b\geq\underline c_R$,  one further has
\[
  |\zeta_i[{\bf y}_j]-\zeta_i[{\bf y}]|
  \leq
  \frac{\norm{c_i[{\bf y}_j]-c_i[{\bf y}]}_{L^2(\Omega)}}
       {\underline c_R|\Omega|^{1/2}}
  \rightarrow0.
\]
Finally, using ${\bf c}[{\bf y}]
  =g\bigl(
    {\bf y}+\boldsymbol{\zeta}[{\bf y}]
    -{\bf z}\phi[{\bf y}]
  \bigr)$
and \eqref{eq:gderivative} in Lemma \ref{lem:localinverse}, we get
\[
\begin{aligned}
  \norm{{\bf c}[{\bf y}_j]-{\bf c}[{\bf y}]}_{L^\infty(\Omega)}
  \leq C\bigl(
    \norm{{\bf y}_j-{\bf y}}_{L^\infty(\Omega)}
    +|\boldsymbol{\zeta}[{\bf y}_j]
      -\boldsymbol{\zeta}[{\bf y}]|+\norm{\phi[{\bf y}_j]-\phi[{\bf y}]}_{L^\infty(\Omega)}
  \bigr)
  \rightarrow0,
\end{aligned}
\]
which yields the continuity of the reconstruction.

\emph{Step 6. Spatial regularity.}
If ${\bf y}\in H^\ell(\Omega)$ with $\ell>1+\frac{d}{2}$, then by using Sobolev embedding, ${\bf y}\in W^{1,\infty}(\Omega)$.  Invoking \eqref{eq:reconstructionbounds}, we further have
$c_p[{\bf y}]-c_n[{\bf y}]\in L^\infty(\Omega)$.
Since $-\Delta\phi[{\bf y}]=c_p[{\bf y}]-c_n[{\bf y}],$ the standard 
elliptic regularity yields, for every $d<r<\infty$,
\[
  \norm{\phi[{\bf y}]}_{W^{2,r}(\Omega)}
  \leq C_r
  \norm{c_p[{\bf y}]-c_n[{\bf y}]}_{L^r(\Omega)}
  \leq C_{r,R},
\]
where $C_r,C_{r,R}>0$ are constants.  The embedding $W^{2,r}(\Omega)\hookrightarrow
W^{1,\infty}(\Omega)$ therefore gives us
$\nabla\phi[{\bf y}]\in L^\infty(\Omega;\mathbb R^d)$. By using the chain rule and \eqref{eq:gderivative}, one has
\[
 \nabla {\bf c}=\bigl[D^2h({\bf c})\bigr]^{-1} (\nabla {\bf y}-{\bf z}(\nabla\phi)^{\top}),\qquad
 \norm{\nabla c}_{L^\infty(\Omega)}\leq\overline{c}_R
                      (\norm{\nabla y}_{L^\infty(\Omega)}+|{\bf z}|\norm{\nabla\phi}_{L^{\infty}(\Omega)}).
\]
\end{proof}

\subsection{A coupled time-discrete scheme}\label{couple-discrete-sec}
Fix an integer $\ell>1+d/2$. For the regularization in entropy variables, set
\[
  Z=\left\{
    w\in H^\ell(\Omega)^2:
    \int_\Omega w_i\,\dx=0,\quad i=n,p
  \right\}.
\]
The regularizing bilinear form is
\begin{align}\label{regular-form-b}
b_\ell({\bf w},{\boldsymbol {\eta}})
  :=\sum_{i=n,p}\sum_{|\beta|\leq\ell}
    \int_\Omega D^\beta w_iD^\beta\boldsymbol{\eta}_i\,\dx.
    \end{align}
 Fix a time step $\tau>0$ and a regularization parameter
$\varepsilon>0$, and take
 ${\bf c}^0={\bf c}_0,$ ${\bf u}^0=P_N{\bf u}_0.$
Given $({\bf c}^{k-1},{\bf u}^{k-1})$, we seek
$({\bf w}^k,{\bf u}^k)\in Z\times V_N$ with $V_N$ given in \eqref{VN-approximation-space}.   Then, one obtains ${\bf c}^k={ c}[{\bf w}^k]$ by using \eqref{eq:reconstructmin} and $\phi^k=\phi[{\bf w}^k]=\invN(c_p^k-c_n^k)$ with $\int_\Omega\phi^k\,\dx=0.$
In addition, we introduce the following regularized discretization
of \eqref{eq:nondimensional-system}:
\begin{align}
 &\frac1\tau({\bf c}^k-{\bf c}^{k-1},\boldsymbol{\eta})
  +(M({\bf c}^k)\nabla {\bf w}^k,\nabla\boldsymbol{\eta})
  -\sum_{i=n,p}(c_i^k{\bf u}^k,\nabla\eta_i)
  +\varepsilon b_\ell({\bf w}^k,\boldsymbol{\eta})=0,
  \qquad \boldsymbol{\eta}\in Z,
  \label{eq:discreteions}\\
  &\frac1\tau({\bf u}^k-{\bf u}^{k-1},{\bf v})
  +\nu(\nabla {\bf u}^k,\nabla {\bf v})
  +b_{\rm NS}({\bf u}^k,{\bf u}^k,{\bf v})
  +\sum_{i=n,p}(c_i^k\nabla {\bf w}_i^k,{\bf v})=0,
  \qquad {\bf v}\in V_N,
  \label{eq:discretems}
 \end{align}
where $(\cdot,\cdot)$ denotes the $L^2$ inner product. Since $\boldsymbol\mu^k=\mathbf w^k+\boldsymbol\zeta^k$,
where $\boldsymbol\zeta^k$ is constant in space, we have
$\nabla\mathbf w^k=\nabla\boldsymbol\mu^k$.
By \eqref{nondimensional-pontential},
\[
 \sum_{i=n,p}c_i^k\nabla w_i^k
 =\nabla\!\left(c_n^k+c_p^k+\tfrac12(\mathbf c^k)^\top F\mathbf c^k\right)
  +(c_p^k-c_n^k)\nabla\phi^k .
\]
The gradient pairs to zero with every $\mathbf v\in V_N$.
Thus \eqref{eq:discretems} has the electrostatic force
$-(c_p^k-c_n^k)\nabla\phi^k$.

\begin{lemma} 
\label{lem:step}
Let $N\in\mathbb N$ and $\tau,\varepsilon>0$.
Assume that
\[
  {\bf c}^{k-1}\in\mathcal K_m,\qquad
  {\bf u}^{k-1}\in V_N,\qquad
  \E({\bf c}^{k-1},{\bf u}^{k-1})<\infty,
\]
where $\E$ is given in \eqref{eq:scaled-total-energy}.  Then there exists a pair
$({\bf w}^k,{\bf u}^k)\in Z\times V_N$ such that, with
${\bf c}^k=c[{\bf w}^k]$ and $\phi^k=\phi[{\bf w}^k]$, \eqref{eq:discreteions} and \eqref{eq:discretems} hold for every
$\boldsymbol{\eta}\in Z$ and ${\bf v}\in V_N$.
\end{lemma} 

\begin{proof}
Fix $N$, $\tau$ and $\varepsilon$. We obtain a fixed point of a compact
solution map by Leray--Schauder degree.

\emph{Step 1. The solution map.}
Set $X=L^\infty_0(\Omega)^2\times V_N$.
Given $(\mathbf y,\widehat{\mathbf v})\in X$ and
$\theta\in[0,1]$, let
$\widehat{\mathbf c}=\mathbf c[\mathbf y]$.
We seek $(\mathbf w,\mathbf v)\in Z\times V_N$ satisfying
\begin{align}
  &\theta\bigl(
    M(\widehat{\mathbf c})\nabla\mathbf w,
    \nabla\boldsymbol\eta
  \bigr)
  -\theta\bigl(
    \widehat{\mathbf c}\otimes\mathbf v,
    \nabla\boldsymbol\eta
  \bigr)
  +\varepsilon b_\ell(\mathbf w,\boldsymbol\eta)
  =-\frac{\theta}{\tau}
    (\widehat{\mathbf c}-\mathbf c^{k-1},
     \boldsymbol\eta),
  \qquad \boldsymbol\eta\in Z,
  \label{eq:homotopyions}\\
  &\frac1\tau(\mathbf v,\boldsymbol\xi)
  +\nu(\nabla\mathbf v,\nabla\boldsymbol\xi)
  +\theta b_{\rm NS}
    (\widehat{\mathbf v},\mathbf v,\boldsymbol\xi)
  +\theta\bigl(
    (\nabla\mathbf w)^\top\widehat{\mathbf c},
    \boldsymbol\xi
  \bigr)
  =\frac1\tau(\mathbf u^{k-1},\boldsymbol\xi),
  \qquad \boldsymbol\xi\in V_N.
  \label{eq:homotopyfluid}
\end{align}
The left-hand sides define a bounded bilinear form on $Z\times V_N$,
since $\widehat{\mathbf c}\in L^\infty(\Omega)^2$ and $V_N$ is finite-dimensional.
Testing with $(\mathbf w,\mathbf v)$ cancels the two coupling terms.
The convection term vanishes by antisymmetry of $b_{\rm NS}$.
Since $M(\widehat{\mathbf c})$ is positive semidefinite, the sum of the
left-hand sides at these tests is bounded below by
\begin{equation}\label{eq:homotopycoercivity}
 \varepsilon\|\mathbf w\|_{H^\ell(\Omega)}^2
 +\tau^{-1}\|\mathbf v\|_{L^2(\Omega)}^2
 +\nu\|\nabla\mathbf v\|_{L^2(\Omega)}^2 .
\end{equation}
This gives coercivity on $Z\times V_N$ uniformly in the input and in $\theta$.
The right-hand sides are bounded linear functionals by Cauchy--Schwarz.
The Lax--Milgram theorem therefore gives a unique solution.
We define the solution map $S_\theta(\mathbf y,\widehat{\mathbf v})
  :=(\mathbf w,\mathbf v)$, where 
 $(\mathbf w,\mathbf v)\in X$.

\emph{Step 2. Continuity and compactness.}
Testing \eqref{eq:homotopyions} with $\mathbf w$ and
\eqref{eq:homotopyfluid} with $\mathbf v$, adding the
two identities, and using \eqref{eq:homotopycoercivity},
we obtain
\begin{align}\label{proceding-bound-homo}
  \|\mathbf w\|_{H^\ell(\Omega)}+\|\mathbf v\|_V
  \leq C\left(
    \|\widehat{\mathbf c}\|_{L^2(\Omega)}
    +\|\mathbf c^{k-1}\|_{L^2(\Omega)}
    +\|\mathbf u^{k-1}\|_{L^2(\Omega)}
  \right).
\end{align}
Thanks to Lemma \ref{lem:reconstruct},
$\|\widehat{\mathbf c}\|_{{L^2(\Omega)}}$ is bounded whenever
$\|\mathbf y\|_{L^\infty(\Omega)}$ is bounded.
Hence, $S_\theta$ maps bounded subsets of $X$ into
bounded subsets of $Z\times V_N$, uniformly for
$\theta\in[0,1]$.
Since $H^\ell(\Omega)\Subset L^\infty(\Omega)$ with $\ell>1+\frac{d}{2}$ and $V_N$ is finite-dimensional, these images are
relatively compact in $X$.

By Lemma~\ref{lem:reconstruct}, the coefficients and
right-hand side of
\eqref{eq:homotopyions}--\eqref{eq:homotopyfluid}
depend continuously on
$(\theta,\mathbf y,\widehat{\mathbf v})$.
Together with the uniform coercivity
\eqref{eq:homotopycoercivity}, continuous dependence
in the Lax--Milgram theorem yields continuity of
$S_\theta$.

\emph{Step 3. A uniform bound for fixed points.}
Suppose $(\mathbf w,\mathbf v)
  =S_\theta(\mathbf w,\mathbf v)$ and $\mathbf c=c[\mathbf w]$ is defined by \eqref{eq:reconstructmin}.     By using \eqref{eq:reconstruction},
$\boldsymbol\mu=\mathbf w+\boldsymbol\zeta$,
where $\boldsymbol\zeta$ is constant in space.
Since $\mathbf c$ and $\mathbf c^{k-1}$ have the
same masses, one obtains from  Lemma \ref{lem:bregman} that  
\[
  (\mathbf c-\mathbf c^{k-1},\mathbf w)
  =(\mathbf c-\mathbf c^{k-1},\boldsymbol\mu)
  \geq
 \Hrel(\mathbf c,0)-\Hrel(\mathbf c^{k-1},0),
\]
where $\Hrel$ is defined in \eqref{H-energy}.  Moreover, $\nabla\mathbf w=\nabla\boldsymbol\mu$, so 
\begin{align}\label{dissipation-equality}
(M(\mathbf c)\nabla\mathbf w,\nabla\mathbf w)
=\D(\mathbf c),
\end{align}
where $\D(\mathbf c)$ is given in  \eqref{nondimensional-dissipation-ionic}.

Taking $\boldsymbol{\eta}=\mathbf w$ in \eqref{eq:homotopyions} and $\boldsymbol{\xi}=\mathbf v$ in \eqref{eq:homotopyfluid}, adding the resulting identities, and multiplying by $\tau$, we use the identity
 $(\mathbf v-\mathbf u^{k-1},\mathbf v)
=\frac12\left(
\|\mathbf v\|_{L^2(\Omega)}^2-\|\mathbf u^{k-1}\|_{L^2(\Omega)}^2
+\|\mathbf v-\mathbf u^{k-1}\|_{L^2(\Omega)}^2
\right)$
to obtain
\begin{equation}\label{eq:homotopyenergy}
\begin{aligned}
  &\theta \Hrel(\mathbf c,0)
  +\frac12\|\mathbf v\|_{L^2(\Omega)}^2
  +\theta\tau\D(\mathbf c)
  +\tau\nu\|\nabla\mathbf v\|_{L^2(\Omega)}^2
  +\tau\varepsilon\|\mathbf w\|_{H^\ell(\Omega)}^2\\
  &\qquad\leq
  \theta \Hrel(\mathbf c^{k-1},0)
  +\frac12\|\mathbf u^{k-1}\|_{L^2(\Omega)}^2
  \leq
  \Hrel(\mathbf c^{k-1},\mathbf u^{k-1}),\end{aligned}
\end{equation}
where we used \eqref{dissipation-equality} and $\Hrel$ is given in \eqref{H-energy}.
Consequently, we have each fixed point satisfies
\[
  \|\mathbf w\|_{H^\ell(\Omega)}^2
  \leq
  \frac{\Hrel(\mathbf c^{k-1},\mathbf u^{k-1})}
       {\tau\varepsilon},
  \qquad
  \|\nabla\mathbf v\|_{L^2(\Omega)}^2
  \leq
  \frac{\Hrel(\mathbf c^{k-1},\mathbf u^{k-1})}
       {\tau\nu}.
\]
These bounds hold uniformly for $\theta\in[0,1]$.
By the embedding $Z\times V_N\hookrightarrow X$,
we can choose $R>0$ such that every fixed point
belongs to the open ball $B_R\subset X$.

\emph{Step 4. Existence of a fixed point.}
At $\theta=0$, \eqref{eq:homotopyions} gives $\mathbf w=0$ and
\eqref{eq:homotopyfluid} determines $\mathbf v\in V_N$ uniquely from
\[
 (\mathbf v,\boldsymbol\xi)+\tau\nu(\nabla\mathbf v,\nabla\boldsymbol\xi)
 =(\mathbf u^{k-1},\boldsymbol\xi)
 \qquad\text{for every }\boldsymbol\xi\in V_N.
\]
Both components are independent of the input.
Thus, $S_0$ is a constant map whose value lies in $B_R$
and $\deg(I-S_0,B_R,0)=1.$
By Step 3, $I-S_\theta$ has no zero on $\partial B_R$
for any $\theta\in[0,1]$. The homotopy invariance therefore
indicates $\deg(I-S_1,B_R,0)=1.$
Hence, $S_1$ has a fixed point.
At this fixed point,
$\widehat{\mathbf c}=\mathbf c[\mathbf w]$ and
$\widehat{\mathbf v}=\mathbf v$, so
\eqref{eq:homotopyions}--\eqref{eq:homotopyfluid}
reduce to
\eqref{eq:discreteions}--\eqref{eq:discretems}.
\end{proof}
Testing the ionic and momentum equations together cancels the work of the
electrostatic force and gives the discrete energy inequality.

\begin{lemma}
\label{lem:discreteenergy}
Let $(\mathbf w^k,\mathbf u^k)\in Z\times V_N$
satisfy \eqref{eq:discreteions}--\eqref{eq:discretems},
with $\mathbf c^k=c[\mathbf w^k]$.
Then
\begin{align}
  &\Hrel(\mathbf c^k,\mathbf u^k)
  +\tau\D(\mathbf c^k)
  +\tau\nu\|\nabla\mathbf u^k\|_{L^2(\Omega)}^2
  +\tau\varepsilon b_\ell(\mathbf w^k,\mathbf w^k)
  \notag\\
  &\qquad
  +\frac12\|\mathbf u^k-\mathbf u^{k-1}\|_{L^2(\Omega)}^2
  +\frac{1}{2}
    \|F^{1/2}(\mathbf c^k-\mathbf c^{k-1})\|_{L^2(\Omega)}^2
  \leq
  \Hrel(\mathbf c^{k-1},\mathbf u^{k-1}),
  \label{eq:discreteenergy}
\end{align}
where $\Hrel$ and $\D$ are given in \eqref{H-energy} and \eqref{nondimensional-dissipation-ionic}, respectively.
\end{lemma}

\begin{proof}
By using \eqref{eq:reconstruction} in Lemma \ref{lem:reconstruct}, we have
$\boldsymbol\mu^k=\mathbf w^k+\boldsymbol\zeta^k$,
where $\boldsymbol\zeta^k$ is constant in space.
Thus, $\nabla\mathbf w^k=\nabla\boldsymbol\mu^k$.
Testing \eqref{eq:discreteions} with $\boldsymbol\eta=\mathbf w^k$
and \eqref{eq:discretems} with $\mathbf v=\mathbf u^k$,
then adding the resulting identities and multiplying by $\tau$,
we obtain \[
\begin{aligned}
  (\mathbf c^k-\mathbf c^{k-1},\mathbf w^k)
  +(\mathbf u^k-\mathbf u^{k-1},\mathbf u^k)
  +\tau\D(\mathbf c^k)
  +\tau\nu\|\nabla\mathbf u^k\|_{L^2(\Omega)}^2
  +\tau\varepsilon b_\ell(\mathbf w^k,\mathbf w^k)=0,
\end{aligned}
\]
where we used 
$b_{\rm NS}(\mathbf u^k,\mathbf u^k,\mathbf u^k)=0$. Since both species masses are conserved,
Lemma~\ref{lem:bregman} yields
\[
\begin{aligned}
  (\mathbf c^k-\mathbf c^{k-1},\mathbf w^k)
  =(\mathbf c^k-\mathbf c^{k-1},\boldsymbol\mu^k)
  \geq
  \Hrel(\mathbf c^k,0)-\Hrel(\mathbf c^{k-1},0)
  +\frac{1}{2}
    \|F^{1/2}(\mathbf c^k-\mathbf c^{k-1})\|_{L^2(\Omega)}^2.
\end{aligned}
\]
Combining this inequality with $(\mathbf u^k-\mathbf u^{k-1},\mathbf u^k)
=\frac12\left(
\|\mathbf u^k\|_{L^2(\Omega)}^2
-\|\mathbf u^{k-1}\|_{L^2(\Omega)}^2
+\|\mathbf u^k-\mathbf u^{k-1}\|_{L^2(\Omega)}^2
\right)$
and the definition of $\Hrel$ in \eqref{H-energy}
proves \eqref{eq:discreteenergy}.
\end{proof}
With the discrete solutions constructed at each time step
in Lemma~\ref{lem:step} and the energy dissipation estimate
from Lemma~\ref{lem:discreteenergy}, we define for $t\in((k-1)\tau,k\tau]$,  $\mathbf c_{\tau,N,\varepsilon}(t)=\mathbf c^k,$
$\mathbf u_{\tau,N,\varepsilon}(t)=\mathbf u^k$ and extend $\mathbf A_i,\mathbf B,\mathbf J_i$
as piecewise constant functions of time in the same way. Then, denote their backward differences as
\begin{align}\label{difference-tau-c}
  \partial_\tau\mathbf c(t)
  :=\frac{\mathbf c^k-\mathbf c^{k-1}}{\tau},
  \qquad
  \partial_\tau\mathbf u(t)
  :=\frac{\mathbf u^k-\mathbf u^{k-1}}{\tau}.
\end{align}
Choose sequences $\tau_j\downarrow0$, $
  \varepsilon_j\downarrow0$, $N_j\uparrow\infty$ 
and denote the corresponding approximations by a
subscript $j$.  Fix $T>0$ and let $K_j=\lceil T/\tau_j\rceil$, so  
$T\leq K_j\tau_j<T+\tau_j$.
\begin{lemma}
\label{lem:approximation-estimates}
Let $\ell>1+d/2$ be an integer and let
$(\mathbf c_0,\mathbf u_0)\in\mathcal K_m\times L^2_\sigma(\Omega)$
have finite relative energy. Consider the piecewise constant interpolants from
Lemma~\ref{lem:step} with $\mathbf c_j^0=\mathbf c_0$ and
$\mathbf u_j^0=P_{N_j}\mathbf u_0$. On every $(0,T)$ they satisfy
\eqref{eq:c-reg}--\eqref{eq:force-reg} uniformly in $j$. Their time differences satisfy
\[
  \|\partial_{\tau_j}\mathbf c_j\|_{L^2(0,T;(H^\ell(\Omega)^2)')}
  +\|\partial_{\tau_j}\mathbf u_j\|_{L^{4/3}(0,T;V')}\le C_T.
\]
The regularization term $\varepsilon_j\mathbf R_j$ defined by
\eqref{eq:Rj-def} tends to zero in $L^2(0,T;(H^\ell(\Omega)^2)')$
at rate $O(\sqrt{\varepsilon_j})$.

Along a common subsequence,
\[
\begin{aligned}
  \mathbf c_j&\to\mathbf c &&\text{strongly in }L^q(0,T;L^2(\Omega)^2)
       \quad(1\le q<\infty),\\
  \mathbf u_j&\to\mathbf u &&\text{strongly in }L^2(Q_T;\mathbb R^d),\\
  \phi_j&\to \invN(c_p-c_n) &&\text{strongly in }L^2(0,T;H^2(\Omega)).
\end{aligned}
\]
Moreover, $\mathbf c_j\to\mathbf c$ almost everywhere in $Q_T$ and
$\mathbf c(t)\in\mathcal K_m$ for almost every $t$. The limits inherit
the bounds on $\mathbf c$, $\mathbf u$ and $\phi$ in \eqref{eq:c-reg}--\eqref{eq:flux-reg}.
\end{lemma}

\begin{proof}
By \eqref{eq:projectionbounds}, projection of the initial velocity can
only decrease the initial energy. Summing \eqref{eq:discreteenergy}
up to step $r$, $1\le r\le K_j$, gives
\begin{equation}\label{eq:uniformsums}
\begin{aligned}
  &\Hrel(\mathbf c_j^r,\mathbf u_j^r)
  +\tau_j\sum_{k=1}^r
  \left[
    \D(\mathbf c_j^k)
    +\nu\|\nabla\mathbf u_j^k\|_{L^2(\Omega)}^2
    +\varepsilon_j
      b_\ell(\mathbf w_j^k,\mathbf w_j^k)
  \right]\\
  &\quad
  +\sum_{k=1}^r
  \left[
    \frac12\|\mathbf u_j^k-\mathbf u_j^{k-1}\|_{L^2(\Omega)}^2
    +\frac{1}{2}
      \|F^{1/2}(\mathbf c_j^k-\mathbf c_j^{k-1})\|_{L^2(\Omega)}^2
  \right]
  \leq \Hrel(\mathbf c_0,\mathbf u_0).
\end{aligned}
\end{equation}
Estimate~\eqref{eq:weighted} now gives
$\tau_j\sum_{k=1}^{K_j}I_w(\mathbf c_j^k)\le C_T$ since $K_j\tau_j<T+\tau_j$.
Lemma~\ref{lem:mixed} yields the stated spatial bounds.
In particular, using also $H^2(\Omega)\hookrightarrow W^{1,6}(\Omega)$ for the electric force, we have
\begin{equation}\label{eq:compactness-spatial-bounds}
\begin{aligned}
 &\|\mathbf c_j\|_{L^\infty(0,T;L^2(\Omega))}
 +\|\mathbf c_j\|_{L^{4/3}(0,T;W^{1,4/3}(\Omega))}
 +\|\mathbf u_j\|_{L^\infty(0,T;L^2(\Omega))}
 +\|\mathbf u_j\|_{L^2(0,T;V)}
 \leq C_T,\\
 &\sum_{i=n,p}
 \left(
   \|c_{i,j}\mathbf u_j\|_{L^2(0,T;L^{4/3}(\Omega))}
   +\|\mathbf J_{i,j}\|_{L^2(0,T;L^{4/3}(\Omega))}
 \right)+\|(c_{p,j}-c_{n,j})\nabla\phi_j\|_{L^2(0,T;L^{3/2}(\Omega))}
 \leq C_T.
\end{aligned}
\end{equation}
For $\boldsymbol\eta\in H^\ell(\Omega)^2$, define
\begin{equation}\label{eq:Rj-def}
  \langle\mathbf R_j(t),\boldsymbol\eta\rangle
  =b_\ell(\mathbf w_j(t),\boldsymbol\eta-\overline{\boldsymbol\eta}).
\end{equation}
By \eqref{eq:uniformsums} and Cauchy--Schwarz in \eqref{regular-form-b},
\begin{equation}\label{eps-wj-estimate}
  \varepsilon_j\|\mathbf w_j\|_{L^2(0,T;H^\ell(\Omega)^2)}^2\le C_T,
  \qquad
  \|\mathbf R_j(t)\|_{(H^\ell(\Omega)^2)'}
  \le C\|\mathbf w_j(t)\|_{H^\ell(\Omega)^2}.
\end{equation}
These bounds give the asserted rate. By discrete mass conservation,
testing \eqref{eq:discreteions} with $\psi-\overline\psi$ gives, for every $\psi\in H^\ell(\Omega)$,
\begin{align}\label{psi-c-approximation}
 \langle\partial_{\tau_j}c_{i,j},\psi\rangle
 =
 (c_{i,j}\mathbf u_j+\mathbf J_{i,j},\nabla\psi)
 -\varepsilon_j\langle R_{i,j},\psi\rangle,
\end{align}
where $R_{i,j}$ is the $i$th component of $\mathbf R_j$.
Since $H^\ell(\Omega)\hookrightarrow W^{1,4}(\Omega)$,
\[
\begin{aligned}
 \|\partial_{\tau_j}\mathbf c_j(t)\|_{(H^\ell(\Omega)^2)'}
 \leq{}&
 C\sum_{i=n,p}
 \left(
   \|c_{i,j}(t)\mathbf u_j(t)\|_{L^{4/3}(\Omega)}
   +\|\mathbf J_{i,j}(t)\|_{L^{4/3}(\Omega)}
 \right)
 +\varepsilon_j\|\mathbf R_j(t)\|_{(H^\ell(\Omega)^2)'}.
\end{aligned}
\]
Squaring and integrating in time using \eqref{eq:compactness-spatial-bounds} and \eqref{eps-wj-estimate} gives
\begin{equation}\label{eq:discretetimeder}
 \|\partial_{\tau_j}\mathbf c_j\|_{L^2(0,T;(H^\ell(\Omega)^2)')}
 \leq C_T.
\end{equation}
Combining \eqref{eq:fluid-L4-prelim} with \eqref{eq:uniformsums} gives
\begin{equation}\label{eq:lady}
 \int_0^T\|\mathbf u_j\|_{L^4(\Omega)}^{8/3}\,dt
 \leq
 C\|\mathbf u_j\|_{L^\infty(0,T;L^2(\Omega))}^{2/3}
 \|\nabla\mathbf u_j\|_{L^2(0,T;L^2(\Omega))}^2
 \leq C_T.
\end{equation}
For $\mathbf v\in V$, using $P_{N_j}\mathbf v$ as a test in \eqref{eq:discretems}, we have from the orthogonality of projection that
\begin{align}\label{test-NS-eq-v}
\begin{aligned}
 \langle\partial_{\tau_j}\mathbf u_j,\mathbf v\rangle
 ={}
 -\nu(\nabla\mathbf u_j,\nabla P_{N_j}\mathbf v)
 -b_{\rm NS}(\mathbf u_j,\mathbf u_j,P_{N_j}\mathbf v)-((c_{p,j}-c_{n,j})\nabla\phi_j,P_{N_j}\mathbf v).
\end{aligned}
\end{align}
By \eqref{eq:projectionbounds}, \eqref{eq:ladyexpanded} and
$V\hookrightarrow L^4(\Omega;\mathbb R^d)$,
\[
 \|\partial_{\tau_j}\mathbf u_j\|_{V'}
 \leq C\left(
   \nu\|\nabla\mathbf u_j\|_{L^2(\Omega)}
   +\|\mathbf u_j\|_{L^4(\Omega)}^2
   +\|(c_{p,j}-c_{n,j})\nabla\phi_j\|_{L^{4/3}(\Omega)}
 \right),
\]
so \eqref{eq:compactness-spatial-bounds} and \eqref{eq:lady} yield
\begin{equation}\label{eq:utapprox}
 \|\partial_{\tau_j}\mathbf u_j\|_{L^{4/3}(0,T;V')}
 \leq C_T.
\end{equation}

Noting that 
\[
 W^{1,4/3}(\Omega)^2
 \Subset L^2(\Omega)^2
 \hookrightarrow (H^\ell(\Omega)^2)' ,
 \qquad
 V\Subset L^2_\sigma(\Omega)\hookrightarrow V',
\]
the time-discrete compactness criterion \cite{Simon,DreherJungel} applied to
\eqref{eq:compactness-spatial-bounds}, \eqref{eq:discretetimeder} and \eqref{eq:utapprox}
gives strong convergence in $L^{4/3}(0,T;L^2(\Omega)^2)$ for $\mathbf c_j$
and in $L^2(Q_T;\mathbb R^d)$ for $\mathbf u_j$ along a common subsequence.
The uniform $L^\infty(0,T;L^2(\Omega)^2)$ bound and interpolation on $(0,T)$
give convergence of $\mathbf c_j$ in $L^q(0,T;L^2(\Omega)^2)$ for every finite $q\ge1$.
Thus, after a further extraction,
\begin{align}
 \mathbf c_j&\rightarrow\mathbf c
 &&\text{strongly in }L^2(Q_T)^2
   \text{ and almost everywhere in }Q_T,
 \label{eq:ccompact}\\
 \mathbf u_j&\rightarrow\mathbf u
 &&\text{strongly in }L^2(Q_T;\R^d).
 \label{eq:ucompact}
\end{align}
We also have 
\begin{align}\label{uj-weak-converge-u-l2}
\mathbf u_j\rightharpoonup\mathbf u\text{ in }L^2(0,T;V).
\end{align}  Then, weak compactness yields
\begin{align}\label{u-Linf-regularity-weak-from}
 \mathbf c\in
 L^\infty(0,T;L^2(\Omega)^2)
 \cap L^{4/3}(0,T;W^{1,4/3}(\Omega)^2),
 \qquad
 \mathbf u\in
 L^\infty(0,T;L^2_\sigma(\Omega))
 \cap L^2(0,T;V).
\end{align}
Applying the bounded Neumann solution operator $\invN:L^2_0(\Omega)\to H^2(\Omega)$ to the charge difference in \eqref{eq:ccompact} gives
\begin{equation}\label{eq:phicompact}
 \phi_j\rightarrow\phi=\invN(z\cdot\mathbf c)
 \quad\text{strongly in }L^2(0,T;H^2(\Omega)).
\end{equation}
Almost everywhere convergence shown in \eqref{eq:ccompact} then implies $\mathbf c\geq0$.  Discrete mass conservation and \eqref{eq:ccompact} yield
\[
 \left\|\int_\Omega c_i(t)\,dx-\int_\Omega c_{i,0}\,dx\right\|_{L^2(0,T)}
 \leq|\Omega|^{1/2}
      \|c_i-c_{i,j}\|_{L^2(Q_T)}
 \rightarrow0.
\]
Hence, $\int_\Omega c_i(t)\,dx=\int_\Omega c_{i,0}\,dx$ for almost every $t$.
Moreover,  $\|{\bf c}_j-\widetilde{\bf c}_j\|
 _{L^2(0,T;(H^\ell)')}\leq
\tau_j\|\partial_{\tau_j}{\bf c}_j\|
 _{L^2(0,T;(H^\ell)')}$ and
\eqref{eq:discretetimeder} imply, by passage to the limit
in $\partial_t\widetilde{\mathbf c}_j
=\partial_{\tau_j}\mathbf c_j$,
 $\partial_t\mathbf c\in L^2(0,T;(H^\ell(\Omega)^2)').$
Together with ${\bf c}_i\in L^\infty(0,T;L^2(\Omega))$ for $i\in\{n,p\}$ obtained by \eqref{eq:ccompact} and the density
of $H^\ell(\Omega)^2$ in $L^2(\Omega)^2$, one has $\mathbf c\in C_{\rm w}([0,T];L^2(\Omega)^2).$
Thus, for this representative,
\[
 \int_\Omega c_i(t)\,dx=\int_\Omega c_{i,0}\,dx
 \qquad\text{for every }t\in[0,T].
\]
\end{proof}


\subsection{Identification of the limiting equations}\label{identify-subsec}
A concentration may vanish in the limit, so $\nabla\log c_i$ need
not be defined. We identify the ionic flux through the weighted fields
$\mathbf A_i$ and $\mathbf B$ from Definition~\ref{def:weak}.

\begin{lemma} 
\label{lem:flux-identification}
Assume all conditions in Lemma \ref{lem:approximation-estimates} hold.  Let $({\bf c}_j,{\bf u}_j,\phi_j)$ be
the subsequence obtained in
\eqref{eq:ccompact}--\eqref{eq:phicompact}.
  Then there exists a   subsequence such that
\begin{equation}\label{eq:ABweak}
  {\bf A}_{i,j}\rightharpoonup {\bf A}_i,\qquad
  {\bf B}_j\rightharpoonup {\bf B}
  \quad\text{in }L^2(Q_T;\R^d),
  \qquad i=n,p,
\end{equation}
and
\begin{align}\label{eq:Jweak1}
  {\bf J}_{i,j}\rightharpoonup {\bf J}_i
  \quad\text{in }L^2(0,T;L^{4/3}(\Omega)).
\end{align}
Moreover, the limits satisfy
\begin{align}\label{limiting-identify}
\begin{aligned}
  {\bf A}_i
  &=
  2\nabla\chi_i({\bf c})
  +\chi_i({\bf c})
  \bigl(
    \nabla\log\omega({\bf c})
    +(F\nabla {\bf c})_i+z_i\nabla\phi
  \bigr),\\
  {\bf B}
  &=
  \sqrt{\frac{bc_n}{\delta}}\,{\bf A}_n
  +\sqrt{\frac{ac_p}{\delta}}\,{\bf A}_p,\\
  {\bf J}_i
  &=-\chi_i({\bf c}){\bf A}_i-\beta_i({\bf c}){\bf B},~~~~~~z_n=-1,~z_p=1 
\end{aligned}
\end{align}
almost everywhere in $Q_T$. In particular, ${\bf A}_i={\bf J}_i=0$ {almost everywhere on } $\{c_i=0\}.$
Also,
\begin{align}\label{assert-lower-semicontinuity}
  \int_0^T I_w({\bf c}(t))\,dt
  \leq
  \liminf_{j\to\infty}
  \int_0^T I_w({\bf c}_j(t))\,dt
  \leq C_T,
\end{align}
where $I_w$ is given in \eqref{eq:weighted} and $C_T>0$ is a constant.
\end{lemma}

\begin{proof}
The bounds \eqref{eq:chi-reg} and \eqref{eq:flux-reg} give
\eqref{eq:ABweak} and \eqref{eq:Jweak1} by weak compactness.
In the formula for $\mathbf A_i$, the coefficients $\chi_i$ and
$\chi_i/\omega$ are bounded and converge almost everywhere. Against
smooth tests they therefore pass to the limit in products with the
weak concentration gradients. For $\mathbf B$, use
$|\sqrt{s}-\sqrt{t}|^2\le|s-t|$: strong $L^2(Q_T)$ convergence of
$c_{i,j}$ gives strong $L^4(Q_T)$ convergence of $\sqrt{c_{i,j}}$,
which pairs with weak $L^2(Q_T)$ convergence of $\mathbf A_{i,j}$.
Continuity and $\beta_i(\mathbf c_j)^4\le C c_{i,j}^2$ give strong
$L^4$ convergence of $\beta_i(\mathbf c_j)$ by uniform integrability.
Its product with $\mathbf B_j$ therefore converges weakly in $L^{4/3}$.
This identifies \eqref{limiting-identify} as in \cite[Section~3.5]{HKLL}.
On $\{c_i=0\}$ the factors $\chi_i,\beta_i$ vanish, as does
$\nabla\chi_i$ almost everywhere. Hence $\mathbf A_i=\mathbf J_i=0$ there.

Estimates \eqref{eq:weighted} and \eqref{eq:uniformsums} give
$$\left\|
\frac{\nabla{\bf c}_j}{{1+c_{n,j}+c_{p,j}}}
\right\|_{L^2(Q_T;\mathbb R^{2\times d})}^2
=
\int_{Q_T}
\frac{|\nabla c_{n,j}|^2+|\nabla c_{p,j}|^2}
     {(1+c_{n,j}+c_{p,j})^2}\,dx\,dt
\leq C_T$$
for  constant $C_T>0.$ Together with \eqref{eq:ccompact} and dominated convergence,
this yields
\[
\frac{\nabla{\bf c}_j}{{1+c_{n,j}+c_{p,j}}}
\rightharpoonup
\frac{\nabla{\bf c}}{{1+c_n+c_p}}
\quad\text{weakly in }L^2(Q_T;\mathbb R^{2\times d}).
\]
Moreover, \eqref{assert-lower-semicontinuity} follows from the weak semicontinuity
of the squared $L^2$ norm.
\end{proof}

\begin{lemma}
\label{lem:coupling-limits}
Assume all conditions in Lemma \ref{lem:flux-identification} hold. Then, we have
\begin{equation}\label{eq:products}
  c_{i,j}{\bf u}_j\to c_i {\bf u},\qquad
  {\bf u}_j\otimes {\bf u}_j\to {\bf u}\otimes {\bf u}
  \quad\text{strongly in }L^1(Q_T).
\end{equation}
Moreover,
\[
  {\bf u}_j\otimes {\bf u}_j\rightharpoonup {\bf u}\otimes {\bf u}
  \quad\text{in }L^{4/3}(0,T;L^2(\Omega)),
\]
and
\begin{equation}\label{eq:forceconvergence}
(c_{p,j}-c_{n,j})\nabla\phi_j
\to(c_p-c_n)\nabla\phi
\end{equation}
strongly in $L^2(0,T;L^{4/3}(\Omega;\mathbb R^d))$ and hence in $L^2(0,T;V')$.
\end{lemma}

\begin{proof}
Thanks to H\"older's inequality, we use \eqref{eq:ccompact} and \eqref{eq:ucompact} to obtain
\[
\begin{aligned}
  \|c_{i,j}{\bf u}_j-c_i {\bf u}\|_{L^1(Q_T)}
  \leq{}
  \|c_{i,j}-c_i\|_{L^2(Q_T)}
  \|{\bf u}_j\|_{L^2(Q_T)}+
  \|c_i\|_{L^2(Q_T)}
  \|{\bf u}_j-{\bf u}\|_{L^2(Q_T)}
  \rightarrow0.
\end{aligned}
\]
Similarly, we have
\[
\begin{aligned}
  \|{\bf u}_j\otimes {\bf u}_j-{\bf u}\otimes {\bf u}\|_{L^1(Q_T)}
  &\leq
  \bigl(\|{\bf u}_j\|_{L^2(Q_T)}+\|{\bf u}\|_{L^2(Q_T)}\bigr)
  \|{\bf u}_j-{\bf u}\|_{L^2(Q_T)}\rightarrow0,
\end{aligned}
\]
which proves \eqref{eq:products}.
Estimate~\eqref{eq:fluid-L4-prelim} bounds
$\|\mathbf u_j\otimes\mathbf u_j\|_{L^{4/3}(0,T;L^2(\Omega))}
=\|\mathbf u_j\|_{L^{8/3}(0,T;L^4(\Omega))}^2$ uniformly.
The strong $L^1$ convergence identifies the weak limit in this space.
For $d=2$, \eqref{eq:fluid-L4-prelim-2d} gives weak convergence
in $L^2(Q_T)$ as well.

By \eqref{eq:phicompact} and $H^2(\Omega)\hookrightarrow W^{1,4}(\Omega)$,
$\nabla\phi_j\to\nabla\phi$ strongly in $L^2(0,T;L^4(\Omega;\mathbb R^d))$.
Consequently,
\[
\begin{aligned}
  &\|(c_{p,j}-c_{n,j})\nabla\phi_j-(c_{p}-c_{n})\nabla\phi\|
  _{L^1(0,T;L^{4/3}(\Omega))}\\
  &\quad\leq
  \|(c_{p,j}-c_{n,j})-(c_{p}-c_n)\|_{L^2(0,T;L^2(\Omega))}
  \|\nabla\phi_j\|_{L^2(0,T;L^4(\Omega))}\\
  &\qquad+
  \|c_{p}-c_{n}\|_{L^2(0,T;L^2(\Omega))}
  \|\nabla\phi_j-\nabla\phi\|_{L^2(0,T;L^4(\Omega))}
  \rightarrow0.
\end{aligned}
\]
By \eqref{eq:force-reg}, the forces are uniformly bounded in
$L^\infty(0,T;L^{4/3}(\Omega;\mathbb R^d))$. Interpolation in time
upgrades the preceding convergence to strong convergence in
$L^2(0,T;L^{4/3}(\Omega;\mathbb R^d))$.
Since $V\hookrightarrow L^4(\Omega;\mathbb R^d)$, this also gives
convergence in $L^2(0,T;V')$.
\end{proof}

\begin{lemma}
\label{lem:limiting-equations}
Assume all conditions in Lemma \ref{lem:flux-identification} hold.  Then, the limits $({\bf c},{\bf u},\phi)$ satisfy
\eqref{eq:weakions} and \eqref{eq:weakns}, together
with the Poisson equation and the prescribed
boundary conditions in the weak sense.
Moreover,
\begin{align}\label{eq:time-continuity}
\begin{aligned}
  c_i&\in
  W^{1,2}(0,T;(W^{1,4}(\Omega))')
  \cap C_{\rm w}([0,T];L^2(\Omega)),\\
  {\bf u}&\in C_{\rm w}([0,T];L^2_\sigma(\Omega)),
  \qquad
  \partial_t{\bf u}\in L^{4/3}(0,T;V'),~i\in\{n,p\},
\end{aligned}
\end{align}
and
\begin{align}\label{IC-c-u-converge}
\begin{aligned}
  c_i(t)&\to c_{i,0}
  &&\text{in }(W^{1,4}(\Omega))',\\
  c_i(t)&\rightharpoonup c_{i,0}
  &&\text{in }L^2(\Omega),\\
  {\bf u}(t)&\rightharpoonup {\bf u}_0
  &&\text{in }L^2_\sigma(\Omega)
\end{aligned}
\qquad\text{as }t\downarrow0.
\end{align}
\end{lemma}

\begin{proof}
Let $\widetilde{\mathbf c}_j$ and $\widetilde{\mathbf u}_j$ interpolate the
nodal values linearly on each time interval. The difference between either
interpolant and its piecewise constant counterpart is bounded pointwise in
time, in the corresponding dual norm, by $\tau_j$ times the norm of its
backward difference. Hence \eqref{eq:discretetimeder} and \eqref{eq:utapprox} give
\begin{equation}\label{eq:affinevanish}
\begin{aligned}
  \widetilde {\bf c}_j-{\bf c}_j&\to0
  &&\text{in }L^2(0,T;(H^\ell(\Omega)^2)'),\\
  \widetilde {\bf u}_j-{\bf u}_j&\to0
  &&\text{in }L^{4/3}(0,T;V').
\end{aligned}
\end{equation}
In addition, the initial values satisfy $\widetilde {\bf c}_j(0)={\bf c}_0$ and $\widetilde {\bf u}_j(0)=P_{N_j}{\bf u}_0.$

Let $\boldsymbol{\eta}=(\eta_n,\eta_p)$ be a smooth test function
with $\boldsymbol{\eta}(T)=0$.
By using the conservation of the discrete species masses, one has $\langle\partial_t\widetilde c_j,\overline{\boldsymbol{\eta}}\rangle=0.$
Testing \eqref{eq:discreteions} with
$\boldsymbol{\eta}-\overline{\boldsymbol{\eta}}$, we therefore obtain
\[
  \langle\partial_t\widetilde c_j,\boldsymbol{\eta}\rangle
  -\sum_{i=n,p}
  (c_{i,j}{\bf u}_j+{\bf J}_{i,j},\nabla\boldsymbol{\eta}_i)
  +\varepsilon_j b_\ell({\bf w}_j,\boldsymbol{\eta}-\overline{\boldsymbol{\eta}})=0,
\]
where ${\bf J}_{i,j}$ is given by  \eqref{eq:Jweak}.
We integrate this identity in time and by parts to obtain
\begin{equation}\label{eq:discreteweakendpoint}
\begin{aligned}
  &-\int_0^T
  \langle\widetilde c_j,\partial_t\boldsymbol{\eta}\rangle\,dt
  -\sum_{i=n,p}\int_{Q_T}
  (c_{i,j}{\bf u}_j+{\bf J}_{i,j})\cdot\nabla\boldsymbol{\eta}_i\,dx\,dt
 +
  \varepsilon_j\int_0^T
  b_\ell(\boldsymbol{w}_j,\boldsymbol{\eta}-\overline{\boldsymbol{\eta}})\,dt
  =({\bf c}_0,\boldsymbol{\eta}(0)),
\end{aligned}
\end{equation}
where $b_\ell$ is given by   \eqref{regular-form-b}.
In light of \eqref{eps-wj-estimate}, we obtain 
\[
\begin{aligned}
  \left|
  \varepsilon_j\int_0^T
  b_\ell({\bf w}_j,\boldsymbol{\eta}-\overline{\boldsymbol{\eta}})\,dt
  \right|&\leq
  C\varepsilon_j
  \|{\bf w}_j\|_{L^2(0,T;H^\ell(\Omega)^2)}
  \|\boldsymbol{\eta}\|_{L^2(0,T;H^\ell(\Omega)^2)}\\
  &\leq
  C_T\sqrt{\varepsilon_j}
\|\boldsymbol{\eta}\|_{L^2(0,T;H^\ell(\Omega)^2)}
  \rightarrow0,
\end{aligned}
\]
where $C$ and $C_T$ are positive constants.  By using \eqref{eq:affinevanish},
Lemma \ref{lem:flux-identification} and
Lemma \ref{lem:coupling-limits}, we pass to the
limit in \eqref{eq:discreteweakendpoint} and
 choose only one nonzero component of $\eta$ to get
\[
  -\int_{Q_T}c_i\partial_t\eta_i\,dx\,dt
  -\int_{Q_T}(c_i {\bf u}+{\bf J}_i)\cdot\nabla\eta_i\,dx\,dt
  =(c_{i,0},\eta_i(0)),
\]
which is \eqref{eq:weakions}.

For the equation of velocity, fix $M\in\mathbb N$ and let $\mathbf v\in C^\infty([0,T];V_M)$ satisfy ${\bf v}(T)=0$.
For all sufficiently large $j$, ${\bf v}(t)\in V_{N_j}$.
Integrating \eqref{eq:discretems} in time gives us
\begin{align}\label{eqNS-approximation}
\begin{aligned}
  &-\int_0^T
  \langle\widetilde {\bf u}_j,\partial_t {\bf v}\rangle\,dt
  -\int_{Q_T}({\bf u}_j\otimes {\bf u}_j):\nabla {\bf v}\,dx\,dt
  +\nu\int_{Q_T}\nabla {\bf u}_j:\nabla {\bf v}\,dx\,dt\\
  &\qquad=
  -\int_{Q_T}(c_{p,j}-c_{n,j}
)\nabla\phi_j\cdot {\bf v}\,dx\,dt
  +(P_{N_j}{\bf u}_0,{\bf v}(0)).
\end{aligned}
\end{align}
In \eqref{eqNS-approximation}, the time derivative term converges by
\eqref{eq:affinevanish} and \eqref{eq:ucompact}. Weak convergence of
$\nabla\mathbf u_j$ gives the viscous term. Lemma~\ref{lem:coupling-limits}
identifies the convection and electric force, and
$P_{N_j}\mathbf u_0\to\mathbf u_0$ in $L^2_\sigma(\Omega)$.
The limiting momentum equation thus holds for smooth tests in finite Stokes spans.
Since $|({\bf u}\otimes {\bf u},\nabla {\bf v})|
  \leq\|{\bf u}\|_{L^4(\Omega)}^2\|\nabla v\|_{L^2(\Omega)}$ and
  $\|u\|_{L^4(\Omega)}^2\in L^{4/3}(0,T)$ as shown in Lemma \ref{lem:approximation-estimates}, we have  $\nabla\cdot(\mathbf u\otimes\mathbf u)
  \in L^{4/3}(0,T;V')$. In addition, we obtained that $ -\nu\Delta\mathbf u
  \in L^2(0,T;V')$ and $-(c_p-c_n)\nabla\phi
  \in L^2(0,T;V')$ shown in Lemma \ref{lem:coupling-limits}.  Density of the finite Stokes spans in $V$
therefore extends the identity to the required
smooth solenoidal tests. This proves
\eqref{eq:weakns} holds.

The velocity boundary condition is retained through
${\bf u}\in L^2(0,T;V)$.
Since $\phi=\invN(c_p-c_n)$, the Poisson equation satisfies
\[
  \int_\Omega\nabla\phi(t)\cdot\nabla\psi\,dx
  =
  \int_\Omega(c_p(t)-c_n(t))\psi\,dx
\]
for every $\psi\in H^1(\Omega)$ and almost every $t$.

It remains to verify the time traces.
In light of \eqref{eq:weakions}, we have
\begin{equation}\label{eq:tracederivative}
  \|\partial_t c_i(t)\|_{(W^{1,4}(\Omega))'}
  \leq
  \|c_i(t){\bf u}(t)+{\bf J}_i(t)\|_{L^{4/3}(\Omega)},
\end{equation}
where the right-hand side belongs to $L^2(0,T)$ by using
\eqref{eq:flux-reg}. Since
$c_i\in L^\infty(0,T;L^2(\Omega))$, we obtain from \eqref{eq:tracederivative} that
\begin{align}\label{contniuous-embeddingci-new}
  c_i\in W^{1,2}(0,T;(W^{1,4}(\Omega))')
  \subset C([0,T];(W^{1,4}(\Omega))').
\end{align}
Moreover, for $\psi\in L^2(\Omega)$, choose
$\psi_m\in W^{1,4}(\Omega)$ with
$\psi_m\to\psi$ in $L^2(\Omega)$.
Then
\[
\begin{aligned}
  |(c_i(t)-c_i(s),\psi)|
  \leq{}
  2C_T\|\psi-\psi_m\|_{L^2(\Omega)}
  +
  |\langle c_i(t)-c_i(s),\psi_m\rangle|,
\end{aligned}
\]
where $C_T>0$ is a constant.  For fixed $m$, the second term in the right hand side tends to zero as
$t\to s$. Letting $m\to\infty$ proves $c_i\in C_{\rm w}([0,T];L^2(\Omega)).$  Similarly, in light of \eqref{eq:weakns}, one takes the limit to obtain
$\partial_t{\bf u}\in L^{4/3}(0,T;V')$.
Together with ${\bf u}\in L^\infty(0,T;L_{\sigma}^2(\Omega))$ shown in \eqref{u-Linf-regularity-weak-from}, the embedding theorem yields 
\begin{align}\label{conituity-u-before3p3}
{\bf u}\in C_{\rm w}([0,T];L^2_\sigma(\Omega)).
\end{align}
Finally, noting that the equations on $(0,T)$ hold in the weak sense,  we apply the integration by parts in time in
\eqref{eq:weakions} and \eqref{eq:weakns} to obtain
\[
  c_i(0)=c_{i,0}
  \quad\text{in }(W^{1,4}(\Omega))',
  \qquad
 {\bf u}(0)={\bf u}_0
  \quad\text{in }L^2_\sigma(\Omega),
\]
which implies the  asserted convergence of the initial data shown in \eqref{IC-c-u-converge} by using
continuity properties \eqref{contniuous-embeddingci-new} and \eqref{conituity-u-before3p3}.
\end{proof}

\subsection{Energy convergence and initial traces}
The decay argument in Section~\ref{sect4} requires the energy inequality
starting at almost every positive time. Strong convergence of the energy
at those times allows passage to the limit at both endpoints.

\begin{lemma}
\label{lem:energy-convergence}
Let $({\bf c},{\bf u},\phi)$ be the limit constructed in Lemma \ref{lem:flux-identification}.  For every $T>0$, there exist a  subsequence
and a set $S\subset(0,T)$ with
$|(0,T)\setminus S|=0$ such that for $t\in S$,
\[
  \Hrel({\bf c}_j(t),{\bf u}_j(t))
  \rightarrow
  \Hrel({\bf c}(t),{\bf u}(t)),
\]
where $\Hrel$ is given in \eqref{H-energy}.
The limiting solution satisfies
\begin{equation}\label{eq:energyalltimes}
\begin{aligned}
  &\Hrel({\bf c}(t),{\bf u}(t))
  +\int_0^t
  \Bigl(
    \|{\bf A}_n(r)\|_{L^2(\Omega)}^2+\|{\bf A}_p(r)\|_{L^2(\Omega)}^2
    +\|{\bf B}(r)\|_{L^2(\Omega)}^2+\nu\|\nabla {\bf u}(r)\|_{L^2(\Omega)}^2
  \Bigr)\,dr\\
  &\qquad\leq\Hrel({\bf c}_0,{\bf u}_0),
  \qquad 0\leq t\leq T.
\end{aligned}
\end{equation}
Moreover,
\begin{equation}\label{eq:twotimeenergy}
\begin{aligned}
  &\Hrel({\bf c}(t),{\bf u}(t))
  +\int_s^t
  \Bigl(
    \|{\bf A}_n(r)\|_{L^2(\Omega)}^2+\|{\bf A}_p(r)\|_{L^2(\Omega)}^2
    +\|{\bf B}(r)\|_{L^2(\Omega)}^2+\nu\|\nabla {\bf u}(r)\|_{L^2(\Omega)}^2
  \Bigr)\,dr\\
  &\qquad\leq\Hrel({\bf c}(s),{\bf u}(s)),
  \qquad
  s\in S\cup\{0\},\quad s\leq t\leq T,
\end{aligned}
\end{equation}
 where ${\bf A}_n$, ${\bf A}_p$ and ${\bf B}$ are given in \eqref{eq:A} and \eqref{eq:B}, respectively.
\end{lemma}

\begin{proof}
Fix $T>0$.
By using \eqref{eq:ccompact} and \eqref{eq:ucompact} in Lemma \ref{lem:approximation-estimates},
after extracting a further subsequence, we find there exists
a measurable set $S\subset(0,T)$ with
$|(0,T)\setminus S|=0$ such that
\begin{equation}\label{eq:goodtimes}
  {\bf c}_j(t)\to{\bf c}(t),\qquad
  {\bf u}_j(t)\to{\bf u}(t)
  \quad\text{strongly in }L^2,
  \qquad t\in S.
\end{equation}
 It then follows from \eqref{eq:goodtimes} and
Lemma~\ref{lem:energycontinuity} that the  $\Hrel({\bf c}_j(t),{\bf u}_j(t))
  \rightarrow
  \Hrel({\bf c}(t),{\bf u}(t))$ on $S$.
Since $P_{N_j}\mathbf u_0\to\mathbf u_0$ in $L^2_\sigma(\Omega)$,
Lemma~\ref{lem:energycontinuity} also gives
\begin{equation}\label{eq:initialenergyconvergence}
  \Hrel({\bf c}_0,P_{N_j}{\bf u}_0)
  \rightarrow
  \Hrel({\bf c}_0,{\bf u}_0).
\end{equation}

Fix $s\in S\cup\{0\}$ and $t\in S$ with $s<t$.
For $0<\eta<t-s$ and all sufficiently large $j$,
\[
  \tau_j\left\lceil\frac{s}{\tau_j}\right\rceil
  <s+\eta,
  \qquad
  t\leq\tau_j\left\lceil\frac{t}{\tau_j}\right\rceil.
\]
Summing \eqref{eq:discreteenergy} in Lemma \ref{lem:discreteenergy} over
$k=\lceil s/\tau_j\rceil+1,\ldots,\lceil t/\tau_j\rceil$,
 we use \eqref{eq:Dfe} to get
\begin{align}\label{Dissipation-upperbound-lem39}
\begin{aligned}
  \Hrel({\bf c}_j(t),{\bf u}_j(t))
  &+\int_{s+\eta}^t
  \left(
    \sum_{i=n,p}\|{\bf A}_{i,j}(r)\|_{L^2(\Omega)}^2
    +\|{\bf B}_j(r)\|_{L^2(\Omega)}^2
    +\nu\|\nabla{\bf u}_j(r)\|_{L^2(\Omega)}^2
  \right)\,dr\\
  &\qquad\leq\Hrel({\bf c}_j(s),{\bf u}_j(s)),
\end{aligned}
\end{align}
where $\Hrel({\bf c}_j(0),{\bf u}_j(0))=\Hrel({\bf c}_0,P_{N_j}{\bf u}_0)$.  For $s,t\in S$, \eqref{eq:goodtimes} and
Lemma~\ref{lem:energycontinuity} give us
$\Hrel({\bf c}_j(s),{\bf u}_j(s))
\rightarrow \Hrel({\bf c}(s),{\bf u}(s))$, $\Hrel({\bf c}_j(t),{\bf u}_j(t))
  \rightarrow \Hrel({\bf c}(t),{\bf u}(t)).$
Combining this with \eqref{eq:initialenergyconvergence},
\eqref{eq:ABweak} and \eqref{uj-weak-converge-u-l2} in Lemma \ref{lem:approximation-estimates}, we pass to the limit
in $  \int_{s+\eta}^t
  \Bigl(
    \|{\bf A}_{n,j}(r)\|_{L^2(\Omega)}^2
    +\|{\bf A}_{p,j}(r)\|_{L^2(\Omega)}^2
    +\|{\bf B}_j(r)\|_{L^2(\Omega)}^2
    +\nu\|\nabla{\bf u}_j(r)\|_{L^2(\Omega)}^2
  \Bigr)\,dr$ and
let $\eta\downarrow0$ to prove
\eqref{eq:twotimeenergy} for
$s\in S\cup\{0\}$ and $t\in S$ with $s<t$.  To extend the inequality to every terminal time,
fix $s\in S\cup\{0\}$ and $t\in(s,T]$, and choose
$t_\ell\in S\cap(s,T)$ with $t_\ell\to t$. By using
\eqref{eq:time-continuity} in
Lemma~\ref{lem:limiting-equations} and the weak lower
semicontinuity  of $\Hrel$ obtained in Lemma~\ref{lem:relative-energy}, we obtain
\[
  \Hrel({\bf c}(t),{\bf u}(t))
  \leq\liminf_{\ell\to\infty}
  \Hrel({\bf c}(t_\ell),{\bf u}(t_\ell)).
\]
Thanks to \eqref{Dissipation-upperbound-lem39}, we have  $\|{\bf A}_n(\cdot)\|_{L^2(\Omega)}^2
  +\|{\bf A}_p(\cdot)\|_{L^2(\Omega)}^2
  +\|{\bf B}(\cdot)\|_{L^2(\Omega)}^2
  +\nu\|\nabla{\bf u}(\cdot)\|_{L^2(\Omega)}^2
  \in L^1(0,T),$
the absolute continuity of the integral implies that,
as $t_\ell\to t$,
\[
\begin{aligned}
  &\int_s^{t_\ell}
  \left(
    \|{\bf A}_n(r)\|_{L^2(\Omega)}^2
    +\|{\bf A}_p(r)\|_{L^2(\Omega)}^2
    +\|{\bf B}(r)\|_{L^2(\Omega)}^2
    +\nu\|\nabla{\bf u}(r)\|_{L^2(\Omega)}^2
  \right)\,dr\\
  &\qquad\rightarrow
  \int_s^t
  \left(
    \|{\bf A}_n(r)\|_{L^2(\Omega)}^2
    +\|{\bf A}_p(r)\|_{L^2(\Omega)}^2
    +\|{\bf B}(r)\|_{L^2(\Omega)}^2
    +\nu\|\nabla{\bf u}(r)\|_{L^2(\Omega)}^2
  \right)\,dr.
\end{aligned}
\]
Combining this convergence with $\Hrel({\bf c}(t),{\bf u}(t))
  \leq
  \liminf_{\ell\to\infty}
  \Hrel({\bf c}(t_\ell),{\bf u}(t_\ell)),$
we pass to the limit in \eqref{eq:twotimeenergy}
with terminal time $t_\ell\in S$ and obtain the same
inequality at $t$.
Thus, \eqref{eq:twotimeenergy} holds for every
$s\in S\cup\{0\}$ and every $t\in(s,T]$.
When $t=s$, the integral vanishes and the inequality in \eqref{Dissipation-upperbound-lem39}
holds  as an equality.
Finally, taking $s=0$ and using the prescribed initial
data gives us \eqref{eq:energyalltimes}.
\end{proof}

\begin{lemma}
\label{cor:stronginitial}
The solution of Lemma~\ref{lem:limiting-equations} satisfies
\begin{align}\label{convergence-strong}
\begin{aligned}
  {\bf c}(t)&\to {\bf c}_0
  &&\text{strongly in }L^2(\Omega)^2,\\
 {\bf  u}(t)&\to {\bf u}_0
  &&\text{strongly in }L^2(\Omega;\R^d) ,\\
  \phi(t)&\to\phi_0
  &&\text{strongly in }H^2(\Omega),
\end{aligned}
\qquad\text{as }t\downarrow0,
\end{align}
where $\phi_0=\invN(z\cdot {\bf c}_0)$.
Moreover, $\Hrel({\bf c}(t),{\bf u}(t))
  \rightarrow
  \Hrel({\bf c}_0,{\bf u}_0)$
 \text{as }$t\downarrow0.$
\end{lemma}

\begin{proof}
Combining \eqref{eq:time-continuity} shown in Lemma \ref{lem:limiting-equations} and weak lower semicontinuity
of $\Hrel$ from
Lemma~\ref{lem:relative-energy}, one has $\Hrel({\bf c}_0,{\bf u}_0)
  \leq
  \liminf_{t\downarrow0}\Hrel({\bf c}(t),{\bf u}(t)).$
On the other hand, \eqref{eq:energyalltimes} implies
\[
  \limsup_{t\downarrow0}\Hrel({\bf c}(t),{\bf u}(t))
  \leq\Hrel({\bf c}_0,{\bf u}_0).
\]
Therefore, $\Hrel({\bf c}(t),{\bf u}(t))
  \rightarrow
  \Hrel({\bf c}_0,{\bf u}_0).$  By using the midpoint inequality \eqref{eq:midpointconvex} in Lemma \ref{lem:relative-energy}, we find
\[
\begin{aligned}
  &\frac{1}{8}\|F^{1/2}({\bf c}(t)-{\bf c}_0)\|_{L^2(\Omega)}^2
  +\frac14\sum_{i=n,p}\int_\Omega
\frac{(c_i(t)-c_{i,0})^2}{c_i(t)+c_{i,0}}\,\dx+\frac18\|{\bf u}(t)-{\bf u}_0\|_{L^2(\Omega)}^2\\
  \leq&
  \frac{\Hrel({\bf c}(t),{\bf u}(t))+\Hrel({\bf c}_0,{\bf u}_0)}2
  -
  \Hrel\left(
    \frac{{\bf c}(t)+{\bf c}_0}{2},
    \frac{{\bf u}(t)+{\bf u}_0}{2}
  \right).
\end{aligned}
\]
Since $\frac{c(t)+c_0}{2}\rightharpoonup c_0$ and $\frac{u(t)+u_0}{2}\rightharpoonup u_0$
  \text{in }$L^2$ as shown in  \eqref{IC-c-u-converge} of Lemma \ref{lem:limiting-equations}, we apply the 
weak lower semicontinuity to get
\[
  \liminf_{t\downarrow0}
  \Hrel\left(
    \frac{{\bf c}(t)+{\bf c}_0}{2},
    \frac{{\bf u}(t)+{\bf u}_0}{2}
  \right)
  \geq\Hrel({\bf c}_0,{\bf u}_0).
\]
Combining this with energy convergence shown in Lemma \ref{lem:energy-convergence} , we have
\[
  \limsup_{t\downarrow0}
  \left(
    \frac{1}{8}\|F^{1/2}({\bf c}(t)-{\bf c}_0)\|_{L^2(\Omega)}^2
    +\frac18\|u(t)-u_0\|_{L^2(\Omega)}^2+\frac14\sum_{i=n,p}\int_\Omega
\frac{(c_i(t)-c_{i,0})^2}{c_i(t)+c_{i,0}}\,\dx
  \right)\leq0.
\]
Moreover, by using Cauchy--Schwarz inequality, we have from \eqref{average-m} that
\[
  \|c_i(t)-c_{i,0}\|_{L^1(\Omega)}^2
  \leq 2m|\Omega|
  \int_\Omega
  \frac{(c_i(t)-c_{i,0})^2}{c_i(t)+c_{i,0}}\,\dx
  \rightarrow0,
  \qquad i=n,p.
\]
For every sequence $t_j\downarrow0$, the functions
$|F^{1/2}\mathbf c(t_j)|^2$ are uniformly integrable.  Due to \eqref{eq:F1}, one has $\alpha|\mathbf c(t_j)|^2
  \leq |F^{1/2}\mathbf c(t_j)|^2$.  Hence, we have  $|\mathbf c(t_j)|^2$  are uniformly integrable.
Together with convergence in measure, Vitali's theorem
therefore yields
\[
  \mathbf c(t_j)\to\mathbf c_0
  \quad\text{strongly in }L^2(\Omega)^2.
\]
Since the sequence was arbitrary, this proves
$\mathbf c(t)\to\mathbf c_0$ strongly in
$L^2(\Omega)^2$ as $t\downarrow0$.
This proves the asserted strong convergence of
${\bf c}(t)$ and ${\bf u}(t)$ shown in \eqref{convergence-strong}.  Finally, $\phi(t)=\invN({\bf z}\cdot {\bf c}(t))$ and the boundedness
of $\invN:L^2_0(\Omega)\to H^2(\Omega)$ imply
\[
  \|\phi(t)-\phi_0\|_{H^2(\Omega)}
  \leq C\|{\bf c}(t)-{\bf c}_0\|_{L^2(\Omega)}
  \rightarrow0
\]
for some constant $C>0.$
\end{proof}

\subsection{Global extension and velocity regularity}
\begin{proof}[Proof of Theorem~\ref{thm:weak}]
Choose $\tau_j,\varepsilon_j\downarrow0$ and $N_j\uparrow\infty$.
By Lemma~\ref{lem:step}, the $j$th approximation exists up to time $j$.
Apply Lemmas~\ref{lem:approximation-estimates}--\ref{lem:energy-convergence}
on $[0,1],[0,2],\ldots$, taking nested subsequences. The diagonal
subsequence defines a weak solution on $[0,\infty)$ since its weakly
continuous representatives agree on overlapping intervals.
Lemma~\ref{cor:stronginitial} gives its strong initial traces.

For each integer $k\ge1$, Lemma~\ref{lem:energy-convergence} on $[0,k]$
gives a null set $N_k\subset(0,k)$ such that \eqref{eq:twotimeenergy}
holds with $T=k$ for every $s\in((0,k)\setminus N_k)\cup\{0\}$.
The union $N=\bigcup_{k\ge1}N_k$ is a null set. For
$s\in((0,\infty)\setminus N)\cup\{0\}$ and $t\ge s$, choose $k>t$; then
\eqref{eq:twotimeenergy} holds on $[s,t]$. Since
$\Hrel(\mathbf c,\mathbf u)=\E(\mathbf c,\mathbf u)-\E(\mathbf c_*,0)$
differs from $\E$ by a constant, this is \eqref{eq:strongenergy}.
Now assume $d=2$. The approximate forces are bounded in $L^2(Q_T;\mathbb R^2)$ by
\eqref{eq:forceL2}. Lemma~\ref{lem:coupling-limits} identifies their
weak limit in this space with $-(c_p-c_n)\nabla\phi$.
Estimate~\eqref{eq:fluid-L4-prelim-2d} puts the convection term in
$L^2(0,T;V')$. Hence the momentum equation gives
$\partial_t\mathbf u\in L^2(0,T;V')$.
Together with $\mathbf u\in L^2(0,T;V)$ this gives
$\mathbf u\in C([0,T];L^2_\sigma(\Omega))$.

Fix $0<s<T$. Choose $t_0\in(0,s)$ such that $\mathbf u(t_0)\in V$.
For the given force and initial velocity $\mathbf u(t_0)$, the
two-dimensional Navier--Stokes regularity theorem
\cite[Chapter~III, Theorem~3.10]{Temam} gives a solution in
$C([t_0,T];V)\cap L^2(t_0,T;D(A))$ with time derivative in
$L^2(t_0,T;L^2_\sigma(\Omega))$.
Uniqueness for this forced equation in the two-dimensional energy class
identifies that solution with $\mathbf u$ on $[t_0,T]$.
Restriction to $[s,T]$ and the Stokes estimate
$D(A)\hookrightarrow H^2(\Omega;\mathbb R^2)$ give the assertion.
When $\mathbf u_0\in V$, the same argument applies directly at $t_0=0$.
\end{proof}

\section{Exponential decay of global weak solutions}\label{sect4}
In this section we prove Theorem~\ref{thm:global-stability}.
The ionic entropy-dissipation estimate and Poincar\'e's inequality for the
velocity give a lower bound on the total dissipation. Combined with the
energy inequality \eqref{eq:strongenergy}, this bound yields exponential decay.

\subsection{A static entropy-dissipation inequality}

\begin{lemma}
\label{lem:gs-entropy-dissipation}
For every $R>0$, there exists $\lambda_R>0$, depending
only on $\Omega,m,a,b,\delta,F$ and $R$, such that
\begin{equation}\label{eq:gs-entropy-dissipation}
\|{\bf A}_n\|_{L^2(\Omega)}^2
+\|{\bf A}_p\|_{L^2(\Omega)}^2
+\|{\bf B}\|_{L^2(\Omega)}^2\geq\lambda_R \Hrel({\bf c},0)
\end{equation}
for every ionic concentration vector
${\bf c}=(c_n,c_p)^\top\in W^{1,4/3}(\Omega)^2$
satisfying \eqref{eq:data} and $\Hrel({\bf c},0)\leq R$, provided that $\chi_n,\chi_p\in W^{1,4/3}(\Omega),$
  ${\bf A}_n,{\bf A}_p,{\bf B}\in L^2(\Omega;\R^d).$
Here $\chi_n$ and $\chi_p$ are defined in
\eqref{eq:weights}, ${\bf A}_n$ and ${\bf A}_p$
in \eqref{eq:A}, and ${\bf B}$ in \eqref{eq:B}.
\end{lemma}

\begin{proof}
We adapt the arguments of
Lemmas~4.2--4.3 and Theorem~4.4 in \cite{HKLL}. Since $F$ is only coercive
on the nonnegative cone, the matrix $L(\mathbf c)$ below need not be bounded.
The collective field $\mathbf B$ is therefore used to control the
concentrations, and the full ionic dissipation enters the argument.

\medskip
\noindent
\textit{Step 1: Energy and gradient estimates.}
By using Lemma~\ref{lem:relative-energy},
\eqref{eq:hm-bounds} and standard elliptic regularity, one has
\begin{align}
  \gamma_*\|{\bf c}-{\bf c}_*\|_{L^2(\Omega)}^2
  &\leq\Hrel({\bf c},0)
  \leq C\|{\bf c}-{\bf c}_*\|_{L^2(\Omega)}^2,
  \label{eq:gs-energy-equivalence}\\
  \|\phi\|_{H^2(\Omega)}
  &\leq C\|{\bf c}-{\bf c}_*\|_{L^2(\Omega)}.
  \label{eq:gs-elliptic}
\end{align}

Set
\begin{equation}\label{eq:gs-g}
  Q=\operatorname{diag}(c_n,c_p),\quad
  K({\bf c})=(I+QF)^{-1},\quad
  L({\bf c})=K({\bf c})Q, \quad
  {\bf g}:=
  \begin{pmatrix}
    \sqrt{c_n\omega/(\delta a)}\,{\bf A}_n^\top\\
    \sqrt{c_p\omega/(\delta b)}\,{\bf A}_p^\top
  \end{pmatrix}.
\end{equation}
By applying the Sobolev chain rule
in \eqref{eq:weights}, together with \eqref{eq:A}, one has
\begin{equation}\label{eq:gs-gradient-identity}
  \nabla{\bf c}
  =K({\bf c}){\bf g}-L({\bf c}){\bf z}(\nabla\phi)^\top.
\end{equation}
This identity holds also on $\{c_i=0\}$, where
$\nabla c_i=0$ almost everywhere.
Since $F$ is positive semidefinite and $f_{nn},f_{pp}\ge\alpha$ by
\eqref{eq:F1}, we have
$\det(I+QF)\ge1+\alpha(c_n+c_p)$. Hence $K$ is continuous on
$[0,\infty)^2$ and its entries are bounded by $\max\{1,|f_{np}|/\alpha\}$.
For $c_n,c_p>0$ we have $L(\mathbf c)=(Q^{-1}+F)^{-1}$ and
$0\le L(\mathbf c)\le Q$. By continuity,
\begin{equation}\label{eq:gs-L-bound}
  |L(\mathbf c)|\le c_n+c_p
  \qquad\text{for all }\mathbf c\in[0,\infty)^2.
\end{equation}
For the rank-one matrix in Remark~\ref{rem:F} and $c_n=c_p=t$, the matrix
$L(\mathbf c)$ has the eigenvalue $t$, so $L$ is unbounded in general.

{ 
Thanks to \eqref{eq:cone-gradient-pointwise},   we apply Cauchy--Schwarz inequality and obtain from \eqref{eq:gs-gradient-identity} that there exist positive constants $C$ and $C_R$ such that
\begin{equation}\label{eq:gs-basic-estimates}
\begin{aligned}
  \left\|
    \frac{\nabla\mathbf c}{1+c_n+c_p}
  \right\|_{L^2(\Omega)}
  &\leq C\left(
    \|\mathbf A_n\|_{L^2(\Omega)}
    +\|\mathbf A_p\|_{L^2(\Omega)}
    +\|\nabla\phi\|_{L^2(\Omega)}
  \right),\\
 \Vert {\bf {g}}\Vert_{L^1(\Omega)}
  &\leq C_R\left(
    \|\mathbf A_n\|_{L^2(\Omega)}
    +\|\mathbf A_p\|_{L^2(\Omega)}
  \right),
\end{aligned}
\end{equation}
where we used $\omega\leq C(1+c_n+c_p)$ shown in \eqref{eq:weights} and \eqref{eq:gs-energy-equivalence}  in the second inequality.  The collective field controls the gradient of the ionic pressure through
\eqref{grad-energy-cn-cp-F}. Applying the argument of
\eqref{upper-lower-bound-cone}--\eqref{eq:mixedholder} at a fixed time,
we find for every $M>0$ a constant $C_{R,M}>0$ such that
\begin{equation}\label{eq:gs-cone-compact-bound}
\begin{aligned}
  &\sum_{i=n,p}
  \left(
    \|c_i\|_{L^4(\Omega)}
    +\|c_i\|_{W^{1,4/3}(\Omega)}
  \right)\leq C_{R,M}
  \quad\text{whenever }
  \|{\bf A}_n\|_{L^2(\Omega)}+\|{\bf A}_p\|_{L^2(\Omega)}
  +\|{\bf B}\|_{L^2(\Omega)}\le M.
\end{aligned}
\end{equation}
 
}


\medskip
\noindent
\textit{Step 2: Vanishing dissipation implies convergence
to equilibrium.}
Let ${\bf c}_j$ satisfy the hypotheses of the lemma and
\[
  \|{\bf A}_{n,j}\|_{L^2(\Omega)}^2
  +\|{\bf A}_{p,j}\|_{L^2(\Omega)}^2
  +\|{\bf B}_j\|_{L^2(\Omega)}^2\to0.
\]
The bound \eqref{eq:gs-cone-compact-bound} and the compact embedding
$W^{1,4/3}(\Omega)\Subset L^2(\Omega)$ give strong convergence of the
concentrations along a subsequence. Neumann regularity gives convergence
of the potentials:
\[
  {\bf c}_j\to{\bf c}
  \quad\text{in }L^2(\Omega)^2
  \text{ and almost everywhere},
  \qquad
  \phi_j\to\phi=\invN(c_p-c_n)
  \quad\text{in }H^2(\Omega).
\]
The limit is nonnegative and retains both masses.
Since ${\bf g}_j\to0$ in $L^1(\Omega)$ by \eqref{eq:gs-basic-estimates}
and $K$ is bounded, $K({\bf c}_j){\bf g}_j\to0$ in $L^1(\Omega)$.
The bound \eqref{eq:gs-L-bound}, the $L^4$ bound in
\eqref{eq:gs-cone-compact-bound} and almost everywhere convergence give
$L({\bf c}_j)\to L({\bf c})$ in $L^2(\Omega)$ by Vitali's theorem.
Together with $\nabla\phi_j\to\nabla\phi$ in $L^4(\Omega)$, this allows
passage to the limit in \eqref{eq:gs-gradient-identity} and yields
\begin{equation}\label{eq:gs-zero-drift}
  \nabla{\bf c}
  =-L({\bf c}){\bf z}(\nabla\phi)^\top.
\end{equation}

We now bootstrap. By \eqref{eq:gs-zero-drift} and \eqref{eq:gs-L-bound},
$|\nabla{\bf c}|\le C|{\bf c}||\nabla\phi|$. Since ${\bf c}\in L^4(\Omega)^2$
and $\nabla\phi\in L^6(\Omega;\R^d)$, we get $\nabla{\bf c}\in L^{12/5}(\Omega)$
and hence ${\bf c}\in L^{12}(\Omega)^2$ for $d\le3$. Neumann regularity
gives $\phi\in W^{2,12}(\Omega)\subset W^{1,\infty}(\Omega)$. Thus
$\nabla{\bf c}\in L^{12}(\Omega)$ and
${\bf c}\in W^{1,12}(\Omega)^2\subset L^\infty(\Omega)^2$.  It then follows from 
 \eqref{eq:gs-zero-drift} that  
\[
  \nabla\left(c_i e^{(F{\bf c})_i+z_i\phi}\right)=0,
  \qquad i=n,p.
\]
Thus,
each $c_i$ is bounded away from zero and
$\mu_i=\log c_i+(F{\bf c})_i+z_i\phi$ is constant.
Hence $(\boldsymbol\mu,{\bf c}-{\bf c}_*)=0$.
Lemma~\ref{lem:bregman} and \eqref{H-energy} now give us 
\[
  0\geq
  \Hrel({\bf c},0)
  +\frac{1}{2}
    \|F^{1/2}({\bf c}-{\bf c}_*)\|_{L^2(\Omega)}^2.
\]
Thus, \eqref{eq:coerciveenergy} implies ${\bf c}={\bf c}_*$.
Relative compactness and uniqueness of this limit imply for the whole sequence, 
\[
  {\bf c}_j\to{\bf c}_*
  \quad\text{in }L^2(\Omega)^2,
  \qquad
  \Hrel({\bf c}_j,0)\to0,
\]
where we used 
Lemma~\ref{lem:energycontinuity}.
 
\medskip
\noindent
\textit{Step 3: Estimate \eqref{eq:gs-entropy-dissipation}.}
Suppose, for contradiction, that there exists a sequence
 such that
\begin{equation}\label{hypothesis-contradiction}
  0<\varepsilon_j^2:=\Hrel({\bf c}_j,0)\leq R,
  \qquad
  \frac{
    \|{\bf A}_{n,j}\|_{L^2(\Omega)}^2
    +\|{\bf A}_{p,j}\|_{L^2(\Omega)}^2
    +\|{\bf B}_j\|_{L^2(\Omega)}^2
  }{\varepsilon_j^2}\to0.
\end{equation}
The conclusion of Step~2 yields ${\bf c}_j\to{\bf c}_*$ and
$\varepsilon_j\to0$.
Define
\[
  {\bf q}_j
  :=\frac{{\bf c}_j-{\bf c}_*}{\varepsilon_j},
  \qquad
  \psi_j:=\frac{\phi_j}{\varepsilon_j}
  =\invN(q_{p,j}-q_{n,j}).
\]
Dividing \eqref{eq:gs-basic-estimates} by
$\varepsilon_j$ and using \eqref{eq:gs-cone-compact-bound},
\eqref{eq:gs-energy-equivalence},
\eqref{eq:gs-elliptic} and
\eqref{hypothesis-contradiction}, we find that
${\bf q}_j$ is bounded in $W^{1,4/3}(\Omega)^2$ and
$g_j/\varepsilon_j\to0$ in $L^1(\Omega)$.  In particular, the gradient bound of ${\bf q}_j$ follows from
\[
  \|\nabla\mathbf q_j\|_{L^{4/3}(\Omega)}
  \leq
  \left\|
    \frac{\nabla\mathbf c_j}
    {\varepsilon_j(1+c_{n,j}+c_{p,j})}
  \right\|_{L^2(\Omega)}
  \|1+c_{n,j}+c_{p,j}\|_{L^4(\Omega)}
  \leq C,
\]
where $C>0$ is a constant.  Consequently, 
\[
  {\bf q}_j\to{\bf q}
  \quad\text{in }L^2(\Omega)^2,
  \qquad
  \psi_j\to\psi=\invN(q_p-q_n)
  \quad\text{in }H^2(\Omega),
\]
and
\[
  {\bf c}_j\to{\bf c}_*
  \quad\text{almost everywhere}.
\]
Moreover, $\int_\Omega q_i\,\dx=0$ for $i=n,p$.  Divide \eqref{eq:gs-gradient-identity} by
$\varepsilon_j$ and pass to the limit as in Step~2.
Since
$L({\bf c}_*)=(m^{-1}I+F)^{-1}$, the vector $(m^{-1}I+F)\mathbf q+\mathbf z\psi$
is constant in space. Testing against the mean-zero field $\mathbf q$ and
using the Poisson equation for $\psi$ gives
\[
  0=
  \int_\Omega{\bf q}^\top(m^{-1}I+F){\bf q}\,\dx
  +\|\nabla\psi\|_{L^2(\Omega)}^2.
\]
Hence ${\bf q}=0$.
However, \eqref{eq:gs-energy-equivalence} implies
\[
  1=\frac{\Hrel({\bf c}_j,0)}{\varepsilon_j^2}
  \leq C\|{\bf q}_j\|_{L^2(\Omega)}^2
  \rightarrow0,
\]
a contradiction. This proves \eqref{eq:gs-entropy-dissipation}.
\end{proof}

\begin{proof}[Proof of Theorem~\ref{thm:global-stability}]
Set $R=\Hrel({\bf c}_0,{\bf u}_0)+1$.
By using \eqref{eq:energyalltimes}, one has
\begin{equation}\label{Hupper-bound-R}
  \Hrel({\bf c}(t),0)
  \leq\Hrel({\bf c}(t),{\bf u}(t))
  \leq\Hrel({\bf c}_0,{\bf u}_0)<R.
  \qquad t\geq0.
\end{equation}
The regularity, nonnegativity,and mass conservation 
established in Theorem~\ref{thm:weak}
ensure that Lemma~\ref{lem:gs-entropy-dissipation}
applies to ${\bf c}(t)$ for almost every $t>0$.  Then, let $\lambda_S>0$ satisfy
\[
  \|\nabla{\bf v}\|_{L^2(\Omega)}^2
  \geq\lambda_S\|{\bf v}\|_{L^2(\Omega)}^2,
  \qquad {\bf v}\in H_0^1(\Omega;\R^d),
\]
and set $\kappa=\min\{\lambda_R,2\nu\lambda_S\}$. We have from
Lemma~\ref{lem:gs-entropy-dissipation}, Poincar\'e's
inequality and \eqref{H-energy} that
\begin{equation}\label{eq:gs-total-dissipation}
\begin{aligned}
  &\|{\bf A}_n(t)\|_{L^2(\Omega)}^2
   +\|{\bf A}_p(t)\|_{L^2(\Omega)}^2
   +\|{\bf B}(t)\|_{L^2(\Omega)}^2
   +\nu\|\nabla{\bf u}(t)\|_{L^2(\Omega)}^2\\
  &\qquad\geq
   \lambda_R\Hrel({\bf c}(t),0)
   +\nu\lambda_S\|{\bf u}(t)\|_{L^2(\Omega)}^2\\
  &\qquad\geq
   \kappa\Hrel({\bf c}(t),{\bf u}(t))
\end{aligned}
\end{equation}
for almost every $t>0$.

Combining \eqref{eq:strongenergy} from Theorem~\ref{thm:weak} with
\eqref{eq:gs-total-dissipation}, we obtain
\begin{equation}\label{eq:gs-integral-decay}
\begin{aligned}
  \Hrel({\bf c}(t),{\bf u}(t))
   +\kappa\int_s^t
     \Hrel({\bf c}(r),{\bf u}(r))\,dr
   \leq\Hrel({\bf c}(s),{\bf u}(s))
\end{aligned}
\end{equation}
for $s=0$ and almost every $s>0$ and for every $t\ge s$.
We now follow the argument  in
Remark~4.1 of \cite{HKLL} to arrive at
\[
  \frac{d}{dt}\Hrel({\bf c}(t),{\bf u}(t))
  +\kappa\Hrel({\bf c}(t),{\bf u}(t))\leq0
\]
in the sense of distributions on $(0,\infty)$.
Thus $e^{\kappa t}\Hrel({\bf c}(t),{\bf u}(t))$
has a nonincreasing representative.
Comparing with times tending to zero and using
\eqref{Hupper-bound-R}, we obtain
\begin{equation}\label{eq:gs-energy-decay}
  \Hrel({\bf c}(t),{\bf u}(t))
  \leq\Hrel({\bf c}_0,{\bf u}_0)e^{-\kappa t}
\end{equation}
for almost every $t>0$.  For any $t>0$, choose $s_j\uparrow t$ such that
\eqref{eq:gs-integral-decay} holds with $s=s_j$ and
\eqref{eq:gs-energy-decay} holds at $s_j$; this is possible since both
hold for almost every time.
Thanks to  \eqref{eq:gs-integral-decay} and the nonnegativity
of $\Hrel$, one finds
\[
  \Hrel({\bf c}(t),{\bf u}(t))
  \leq\Hrel({\bf c}(s_j),{\bf u}(s_j))
  \leq\Hrel({\bf c}_0,{\bf u}_0)e^{-\kappa s_j}.
\]
Letting $s_j\uparrow t$ extends \eqref{eq:gs-energy-decay} to every $t\ge0$.
Coercivity and the Poisson estimate \eqref{eq:gs-elliptic} then give
\[
\begin{aligned}
  \|{\bf c}(t)-{\bf c}_*\|_{L^2(\Omega)}
   +\|{\bf u}(t)\|_{L^2(\Omega)}
   +\|\phi(t)\|_{H^2(\Omega)}\leq
   C\Hrel({\bf c}(t),{\bf u}(t))^{1/2}
   \leq
   C\Hrel({\bf c}_0,{\bf u}_0)^{1/2}e^{-\kappa t/2}.
\end{aligned}
\]
This proves \eqref{eq:global-norm-decay}. Here $C$ depends only on
$\Omega,F,m,\alpha$, and $\kappa$ depends only on
$\Omega,m,a,b,\delta,\nu,F$
and $\Hrel({\bf c}_0,{\bf u}_0)$.
\end{proof}

\section{Strong solutions near equilibrium}
\label{sec:small-strong}
We prove Theorem~\ref{thm:small-strong} by adapting the argument in
Section~4.4 of \cite{HKLL} and estimating the additional fluid coupling terms.
\begin{proof}[Proof of Theorem~\ref{thm:small-strong}]
Set $r>\max\{3,d\}$ and write $\mathbf q=\mathbf c-\mathbf c_*$ and $U=(\mathbf q,\mathbf u)$.
Let $A_r$ be the Dirichlet Stokes operator on $L^r_\sigma(\Omega)$ and set
\[
\begin{aligned}
 X_0&=L^r_0(\Omega)^2\times L^r_\sigma(\Omega),\\
 X_1&=\left\{\mathbf q\in W^{2,r}(\Omega)^2:
       \int_\Omega\mathbf q\,\dx=0,\quad
       \partial_{\boldsymbol\nu}\mathbf q=0\right\}\times D(A_r),
 \qquad X_\gamma=(X_0,X_1)_{1-1/r,r}.
\end{aligned}
\]
Here $L^r_0(\Omega)$ denotes the mean-zero subspace and
$D(A_r)=W^{2,r}(\Omega;\mathbb R^d)\cap W^{1,r}_0(\Omega;\mathbb R^d)
\cap L^r_\sigma(\Omega)$. In dimension one the solenoidal component
is zero. Since $2-2/r>1+1/r$, the trace space $X_\gamma$ imposes the
Neumann condition on $\mathbf q$ and the no-slip condition on $\mathbf u$.
Its norm is equivalent to the sum of the two $W^{2-2/r,r}$ norms.
Also $2-2/r>d/r$, so $X_\gamma$ embeds into $L^\infty$ componentwise.

Write $M_*=M(\mathbf c_*)$ and $H_*=H(\mathbf c_*)=m^{-1}I+F$.
For small $\|U\|_{X_\gamma}$ we have $m/2<c_i<3m/2$.
Both $M(\mathbf c)$ and $H(\mathbf c)$ are then invertible, so
\eqref{eq:bc} is equivalent to $\partial_{\boldsymbol\nu}\mathbf q=0$
for the ionic equations. Using $\phi=\invN(\mathbf z\cdot\mathbf q)$,
we write \eqref{eq:nondimensional-system} as
\begin{equation}\label{eq:small-strong-linearization}
 \partial_t U+\mathcal L U=\mathcal N(U),\qquad
 \mathcal L U=
 \bigl(-M_*H_*\Delta\mathbf q+M_*\mathbf z\mathbf z^\top\mathbf q,
       \,\nu A_r\mathbf u\bigr).
\end{equation}
The ionic component of $\mathcal N(U)$ is given, for $i=n,p$, by
\[
 \nabla\cdot\left(
 \sum_{j=n,p}[M(\mathbf c)H(\mathbf c)-M_*H_*]_{ij}\nabla q_j
 +[(M(\mathbf c)-M_*)\mathbf z]_i\nabla\phi\right)
 -\mathbf u\cdot\nabla q_i .
\]
It has zero spatial mean by the boundary conditions and incompressibility.
The fluid component is the Helmholtz projection of
$-(\mathbf u\cdot\nabla)\mathbf u-(\mathbf z\cdot\mathbf q)\nabla\phi$.

The matrix $M_*H_*$ is similar to the positive definite matrix
$M_*^{1/2}H_*M_*^{1/2}$. A constant change of variables thus reduces
the principal ionic operator to two scalar Neumann Laplacians.
The remaining ionic term is bounded on $X_0$, so the Neumann and
Stokes maximal $L^r$ regularity estimates apply to $\mathcal L$.
On each positive Neumann mode with eigenvalue $\lambda$, the ionic
matrix is similar to
$M_*^{1/2}(\lambda H_*+\mathbf z\mathbf z^\top)M_*^{1/2}$.
Its smallest eigenvalue has a positive lower bound uniform in
$\lambda$ because the constant mode is excluded. The Stokes operator
also has a positive spectral gap whenever the fluid component is present.

Choose $\eta>0$ below the spectral gap of $\mathcal L$. Define
\[
 \mathbb E_\eta=\{U:e^{\eta t}U\in
 W^{1,r}(0,\infty;X_0)\cap L^r(0,\infty;X_1)\}
\]
with the sum of the two norms of $e^{\eta t}U$. Maximal regularity for
$\mathcal L-\eta I$ and the time trace theorem
\cite{PrussSimonett2016,PSZ2009} give
\begin{equation}\label{eq:small-strong-maximal}
 \sup_{t\ge0}e^{\eta t}\|U(t)\|_{X_\gamma}
 \le C\|U\|_{\mathbb E_\eta}
 \le C'\left(\|U(0)\|_{X_\gamma}
 +\|e^{\eta t}(\partial_t U+\mathcal L U)\|_{L^r(0,\infty;X_0)}\right).
\end{equation}

To estimate the differentiated coefficients in the ionic flux, apply the
Gagliardo--Nirenberg inequality on $\Omega$ at each time:
\[
 \|\nabla v\|_{L^{2r}(\Omega)}^2
 \le C\|v\|_{L^\infty(\Omega)}\|v\|_{W^{2,r}(\Omega)}.
\]
For any scalar component $v$ of $U$, integration and the trace estimate yield
\[
 \|e^{\eta t/2}\nabla v\|_{L^{2r}((0,\infty)\times\Omega)}^2
 \le C\sup_{t\ge0}\|v(t)\|_{L^\infty(\Omega)}
       \|e^{\eta t}v\|_{L^r(0,\infty;W^{2,r}(\Omega))}
 \le C\|U\|_{\mathbb E_\eta}^2.
\]
Thus products of two first derivatives belong to the required weighted
$L^r$ space. The coefficient of each $\Delta q_j$ in the nonlinear
remainder is bounded by $C|\mathbf q|$. For the electrostatic terms,
elliptic regularity and $r>d$ give
$\|\nabla\phi(t)\|_{L^\infty}\le C\|\mathbf q(t)\|_{L^r}$,
and $-\Delta\phi=\mathbf z\cdot\mathbf q$.
The transport terms are controlled by placing the undifferentiated
velocity in $L^\infty$.

On the concentration range $[m/2,3m/2]^2$, the coefficients and their
first two derivatives are bounded. Subtracting the fluxes, a difference
of coefficients is bounded by $C|\mathbf q-\widehat{\mathbf q}|$.
A difference of gradient products has one differentiated difference in
each term and is controlled by the preceding $L^{2r}$ estimate.
Together these bounds give, on a sufficiently small ball in $\mathbb E_\eta$,
\begin{equation}\label{eq:small-strong-nonlinear}
 \|e^{\eta t}(\mathcal N(U)-\mathcal N(\widehat U))\|_{L^r(0,\infty;X_0)}
 \le C(\|U\|_{\mathbb E_\eta}+\|\widehat U\|_{\mathbb E_\eta})
      \|U-\widehat U\|_{\mathbb E_\eta}.
\end{equation}
Given $\widehat U$, solve the linear equation with right-hand side
$\mathcal N(\widehat U)$ and the prescribed initial data.
Since $\mathcal N(0)=0$, estimates
\eqref{eq:small-strong-maximal}--\eqref{eq:small-strong-nonlinear}
make this map a contraction on a ball of radius
$2C'\|U(0)\|_{X_\gamma}$ for sufficiently small $\varepsilon_*$.
Its fixed point is a global strong solution. The same difference
estimate on short time intervals gives uniqueness in the class
\eqref{eq:small-strong-regularity}: starting at a common trace,
the difference has zero initial value, the trace norms stay small by
continuity, and the time integrals of the $W^{2,r}$ norms tend to zero
with the interval length. Repeating the argument gives uniqueness on
every finite interval.

The trace bound in \eqref{eq:small-strong-maximal} gives the concentration
bounds and the exponential decay of $\mathbf q$ and $\mathbf u$.
The Neumann estimate
$\|\phi(t)\|_{W^{4-2/r,r}(\Omega)}
 \le C\|\mathbf q(t)\|_{W^{2-2/r,r}(\Omega)}$
then gives \eqref{eq:small-strong-decay}.
\end{proof}

\section{Conclusion}\label{sect-final}
In this paper, we derive a PNP--NS system with interspecies drag and
steric effects using the energetic variational approach. Combining
weighted dissipation estimates with a coupled implicit Euler--Galerkin
scheme that preserves the energy-dissipation structure, we establish the
existence of global finite-energy weak solutions in dimensions
$1\leq d\leq3$ for arbitrary finite-energy data, allowing ionic vacuum.
The electrostatic force is square integrable in space and time, and in
two dimensions the solvent velocity becomes regular at positive times. We further prove by the
entropy method that the constructed solutions converge exponentially to
the constant equilibrium. For
sufficiently small perturbations of the equilibrium, we also prove global
existence and uniqueness of strong solutions, together with exponential
decay.

Several questions remain open, including weak--strong
uniqueness and higher regularity of the coupled system. In two dimensions,
global classical solvability for large smooth positive initial data
requires further concentration estimates. The optimal rates at which the
ionic concentrations approach equilibrium and the solvent velocity tends to
zero remain to be determined. Other directions include the stability of
nonconstant steady states under nonhomogeneous boundary conditions and
extensions to multicomponent electrolytes.
\section*{Acknowledgments} 
This work was mostly carried out when Chun Liu was still working in Illinois Institute of Technology.
Baoli Hao and Chun Liu were partially supported by {DMS-2216926 and DMS-2410742}. The authors used ChatGPT (OpenAI) to assist with language editing and reviewing mathematical proofs.
All content was reviewed and verified by the authors, who take full responsibility for it.







\end{document}